\input amstex 
\documentstyle{amsppt}                         
\magnification=1200
\define\={=}
\define\C{{\Bbb C}}
\define\R{{\Bbb R}}

\define\re{{\operatorname{Re}}}
\define\im{{\operatorname{Im}}}
\define\unim{{\sqrt{-1}}}

\redefine\D{{\Cal D}}
\redefine\L{{\Cal L}}
\define\E{{\Cal E}}
\define\W{{\Cal W}}
\redefine\H{{\Cal H}}
\define\F#1{{\Cal F}^{(#1)}}

\define\g{\goth g}
\define\gl{\goth{gl_n(\C)}}
\redefine\u{\goth u_{n-1}\oplus \R}

\redefine\P{{\Bbb P}}
\define\TM{\tilde TM}
\define\PTM{{\Bbb P}TM}
\define\CP{\C P}

\define\LC{L^\C(M)}
\define\UG{U_g(M)}
\define\UF{U_{F}(M)}
\define\ULSM{U_{\rho}(SM)}

\define\Hom{\operatorname{Hom}}
\define\rank{\operatorname{rank}}
\define\bh{\text{h}}
\define\bH{\text{H}}
\define\gh{\text{\bf H}}

\documentstyle{amsppt}
\TagsOnRight

\topmatter

\title
The Structure Equations\\ of a Complex Finsler Manifold
\endtitle

\author
A. Spiro
\endauthor

\address
\phantom{ }\newline
A. Spiro\newline
Dipartimento di Matematica e Fisica\newline
Universit\`a di Camerino\newline
62032 CAMERINO (Macerata)\newline
ITALY\newline
\phantom{ }\newline
\phantom{ }
\endaddress

\email {spiro\@campus.unicam.it}
\endemail


\leftheadtext{A. Spiro}
\rightheadtext{The Structure Equations\\ of a Complex Finsler Manifold}

\abstract For a strongly pseudo-convex complex Finsler manifold
$M$,
a bundle $\UF$ of adapted unitary frames is canonically defined.  
A non-linear Hermitian connection
on $\UF$, invariant under  local  biholomorphic isometries, is given and
it proved to be unique. By means of such connection, an absolute parallelism 
on $\UF$ is determined and
 a new set
of structure functions which generate all the isometric invariants
of a Finsler metric is obtained.\par
 A pseudo-convex complex Finsler manifolds $M$,
which  admits a totally geodesic complex curve with a given constant holomorphic 
sectional curvature through any point and any direction, is called
{\it E-manifold\/}. Main examples of E-manifolds
are the smoothly bounded, strictly convex 
domains in $\C^n$, endowed with the Kobayashi metric.
A complete characterization  of E-manifolds, using the previously defined
structure functions, is given and a smaller set of generating functions
for the isometric invariants of  
E-manifolds is determined.\par

\endabstract

\subjclass Primary 53B40, 32H15; Secondary 53C60, 53A55, 53B15
\endsubjclass
\keywords  Complex Finsler Metrics, Kobayashi Metrics
\endkeywords
\endtopmatter

\document
\subhead 1. Introduction
\endsubhead
\bigskip
The main purpose of this paper is to give a 
complete set of  invariants, which  characterize a strongly pseudoconvex
complex Finsler metric up to local biholomorphic isometries.
Several properties of these invariants are immediately related with
the intrinsic geometry of the
Kobayashi  metric of the smoothly bounded, strongly convex domains
in $\C^n$.
\par
\medskip
Let $M$ be a complex manifold and  $J$  its  complex structure. 
The well-known {\it infinitesimal Kobayashi  pseudo-distance $k_M$ on $\TM = 
TM \setminus \{\text{\rm zero section}\}$\/} can be defined 
 as follows ([Ko]):
for any 
$x\in M$ and any $0\neq v\in T_xM$, let 
 $\Cal A_v$ 
 the set of all 
$r\in \R^+$ such that there exists a  
holomorphic map $f: \Delta_r\to M$ from 
$\Delta_r = \{ |z| <r\}\subset \C$ into $M$ with $f(0) = x$
and $f_*(\frac{\partial}{\partial \zeta}) \in \C v$. Then 
$$k_M(v) = \inf_{r\in \Cal A(v)} \frac{1}{r}\ .$$
We consider  the following
class of complex manifolds.
\par
\definition{Definition 1.1} A complex manifold $(M,J)$ is a
{\it Lempert manifold\/} if 
\roster 
\item the infinitesimal Kobayashi pseudo-distance $k_M$ is 
a strongly pseudoconvex  Finsler metric, that is: \par
a) it is a
smooth function on $\TM$ with values in $\R^+$; \par
b)  
 $k_M(\lambda v) = |\lambda|k_M(v)$ for any $\lambda \in \C^*$ and $v\in \TM$;
\par
 c) at any point $x\in M$ the hypersurface
$S_x = \{ v\in T_xM\ : k_M(v) = 1\}$
is strongly pseudoconvex  in $T_xM$;
\item for any non-vanishing
 complex vector $w\in T^\C_x M\subset T^\C M$, there exists a
complex curve $\gamma_w: U\subset \C \to M$, such that 
$\gamma_w(0) = x$, $\gamma'_w(0) = w$ and 
$\gamma(U)$  is a totally geodesic submanifold of $M$;
\item the metric, which is induced by $k_M$ on the totally geodesic complex
curve $\gamma_w(U)$, is K\"ahler  and 
with constant holomorphic curvature  equal to $-4$;
\item the (finite) Kobayashi distance $d_M$, determined by $k_M$, 
is complete and the exponential map $\exp: T_xM \to M$ is 
a diffeomorphism for any $x\in M$.
\endroster
\enddefinition
An immediate interest for  Lempert manifolds
comes from the  well-known results of  L. Lempert on the Kobayashi
metric of strongly convex
domains in $\C^n$ ([Le], [Le1]). 
Some of  his results can be  stated as follows.\par
\proclaim{Theorem 1.2}[Le] If $M$ is a smoothly bounded, strongly 
convex domain in $\C^n$, then $M$ is a Lempert manifold.
\endproclaim
Since the convexity of a domain is 
not a biholomorphic invariant property,
Theorem 1.2 motivates  the following  
question:\par
{\it Are there some invariant 
properties of $k_M$ (to be added to  (1) - (4) of Def. 1.1), 
which characterize  the manifolds
that
are biholomorphic to a  smoothly bounded, strongly 
convex  domain in $\C^n$?\/}\par
\medskip
Some encouraging   results  have been  obtained by various authors
(e.g. [Fa], [Pa], [Le2], [BD], [AP], [AP1]).
In particular, we would like to mention the following theorem by J. J. Faran
(see also [Pa]).
\par
\proclaim{Theorem 1.3}[Fa] $(M,J)$ is a Lempert manifold if and only if 
it admits a strongly pseudoconvex 
Finsler metric $F$, which verifies (2), (3) and (4) of Definition
1.1. In this case $F$ coincides with the Kobayashi metric $k_M$.
\endproclaim
Faran's Theorem has been improved by M. Abate and G. Patrizio in [AP]
in the following sense: they proved that if 
$(M,J)$ has a strongly pseudoconvex 
Finsler metric $F$, then it admits a natural non-linear Finsler connection
and if  the corresponding
 torsion and curvature  
verify pointwise
a certain set of conditions, then $(M,J)$ satisfies 
(2), (3) and (4) of 
Def. 1.1
and  hence it  is a Lempert manifold.\par 
\medskip
By Faran's theorem,
our previous question has positive answer if and only if there
exist some conditions, which are necessary and
sufficient for the existence of a biholomorphism between a Finsler manifold
$M$, verifying 
(2), (3) and (4) of Def. 2.1, and a strongly convex domain $D\subset \C^n$. \par
\medskip
The general problem of determining necessary and sufficient
conditions  
for the existence of a  (local) isomorphism between two geometric structures
is  usually called {\it the 
(local) equivalence problem\/} for those structures. 
In this paper we give a new solution to
the local equivalence problem for   
strongly pseudoconvex
 Finsler metrics and, by means of this solution,
 we obtain a new complete set 
of invariant functions which determine the complex Finsler metrics
up to local biholomorphic isometries. We also use these invariants
to give a new characterization of Lempert manifolds. 
\medskip
Here are the contents of the paper. In \S 2  we recall and prove some preliminary 
properties of complex Finsler metrics.\par
In \S 3 we introduce  the concept of {\it adapted unitary frames\/}
of a complex  manifold $(M,J)$ with a strongly pseudoconvex Finsler metric  $F$.
The bundle $U_F(M)$ of  all adapted unitary frames
turns out to be a subbundle of the complex linear frame bundle $\LC$, but 
in general  it is  {\it not\/} a principal subbundle; this is the case
if and only if there exists an Hermitian metric $g$ so that $F(v)  =
\sqrt{g(v,v)}$ for all $0\neq v \in TM$. \par
We also use the following terminology:
  any 
distribution  which is complementary to the vertical 
distribution and of dimension equal to $\dim M$ is named a 
{\it non-linear connection on $U_F(M)$\/}. We say that 
a non-linear connection is {\it of Hermitian type\/} if it is invariant 
w.r.t. the complex structure $\hat J$ of  $\LC$. \par
The main result of \S 3 is the following (Theorem 3.10).\par
\proclaim{Theorem 1.4} Let $(M,J,F)$ be a strongly pseudoconvex
Finsler manifold. Then the unitary frame bundle $\UF$ has a unique
non-linear connection of Hermitian type. This connection
is invariant under any biholomorphic isometry of $(M, J, F)$.
\endproclaim
This non-linear Hermitian connection 
on $\UF$ defines a non-linear covariant derivation for vector fields of $M$
which is invariant under any biholomorphic isometry.\par 
In \S 4 we  show that any fiber $\Bbb V_x = \pi^{-1}(x)$ of $\UF$
is identifiable with the adapted frame bundle of a pseudo-hermitian 
structure on the Finsler sphere $S_x$. Using the 
Webster connection for pseudo-hermitian structures (see [We]), 
we  define  an invariant absolute parallelism on each fiber $\Bbb V_x$, 
i.e. a set of vertical vector fields on $\Bbb V_x$, which  at all points 
span $T_u \Bbb V_x$ and 
which is invariant under the automorphism of the pseudo-hermitian structure
of $S_x$. \par
Using this absolute parallelism on the fibers
 and the non-linear Hermitian connection 
$\H$  of $\UF$, we obtain
an absolute parallelism $\sigma$ on $\UF$ 
which  verifies the following crucial property: 
{\it the (local) biholomorphic isometries of $(M,J,F)$ are in 1-1 correspondence 
with the (local) diffeomorphisms of $U_F(M)$ which preserve $\sigma$\/}.\par
By Kobayashi's theorem on the automorphisms of absolute 
parallelisms ([Ko1]), we immediately obtain the following result (Proposition 
4.6):
\medskip
\proclaim{Theorem 1.5} Let $(M,J)$ be a complex manifold of complex 
dimension $n$ and  $F$   a strongly pseudoconvex 
Finsler metric on $(M,J)$. \par
The group of 
biholomorphic isometries $Iso_F(M,J)$ is a Lie group of dimension less or
equal to
$n^2 + 2n$. Moreover $\dim_\R Iso_F(M,J) = n^2 + 2n$ 
 if and only if  $F$ is equal to
$F(v) = \sqrt{g(v,v)}$ for some 
K\"ahler metric $g$ of constant holomorphic sectional
curvature and $(M, J, g)$ is a simply connected complex space form, 
i.e. $\C P^n$, 
$\C^n$ or the unit ball $B^n \subset \C^n$, endowed with a
Fubini-Study, flat  or Poincar\`e-Bergmann metric, respectively 
\endproclaim
In \S 5  we determine the Lie brackets
of all possible pairs  of vector fields  of  
the absolute parallelism 
$\sigma$ 
of $\UF$. By Cartan-Sternberg theorem the components of these Lie brackets
w.r.t. the vectors of the absolute parallelism
 generate a complete 
set of invariant functions for the Finsler manifold $(M,J,F)$ (see
Proposition 4.6 and Theorem A1). At the end of \S 5,
  we also give the so-called {\it structure equations of
the 
 Finsler manifold of $(M,J,F)$\/}, 
i.e. the equations that are verified by 
 the  1-forms on $U_F(M)$ which are dual to the
vector fields of the absolute parallelism. \par
\medskip
At last, in \S 6, we determine the Euler-Lagrange equations for the 
geodesics of a complex Finsler manifold. 
We  recall the definition of complex geodesics (see [Ve] and [AP])
and we find necessary and sufficient conditions for a complex
Finsler manifold to be of constant holomorphic sectional
curvature and with a complex geodesic through any point
and  any direction. We call such manifolds {\it E-manifolds\/}.
Notice that the Lempert manifolds are  complete E-manifolds
with  holomorphic sectional curvature equal to -4. \par
For the E-manifold, we also  prove that
the torsion and  the curvature can be expressed
in terms of  the other structure functions of the absolute
parallelism on $\UF$ and hence that these 
structure functions are the actual generators for
the invariants of  E-manifolds (see  Theorem 6.9).\par
\medskip
We have to mention that an alternative
 solution to the equivalence
problem  has been  given  by J. J. Faran in [Fa]. 
He  determines another set
of  invariant
functions, by pursuing the steps of 
a general algorithmic 
procedure: it is our personal opinion
that,  by this reason, it is quite cumbersome
to obtain simple geometric 
interpretations  for the invariant functions
introduced by Faran. \par
\medskip
As final remark, we want to point out that our  non-linear 
covariant derivation on complex Finsler 
manifolds 
 is 
strictly related (but different) with the non-linear covariant derivation 
invented by S. Kobayashi in [Ko2].  In a forthcoming paper,
 we will discuss this  relation and 
we will show how to use 
 the components of the curvature  and torsion tensors 
in order to determine the set of invariant  functions given in 
\S 5. \par
\bigskip
\subhead 2. Preliminaries 
\endsubhead
\medskip
\subsubhead 2.1 Notation and basic definitions
\endsubsubhead
\medskip
In all this paper, we   use
 greek letters $\alpha$, $\beta$, etc.
for indices related to holomorphic vectors,  barred greek letters
$\bar\alpha$, $\bar \beta$, etc.  for indices related to the  
conjugated
vectors  and latin indices $i, j, k$, etc.
to denote  real vectors.\par
We  denote by  $\{\epsilon_0, \epsilon_1, \dots, 
\epsilon_{2n-1}\}$ 
the standard real basis of $V= \R^{2n} = \C^n$; $J_o$
is the  complex structure of $\C^n$. The standard basis is ordered
so that 
$J_o(\epsilon_{2i}) = \epsilon_{2i + 1}$ for any $i = 0, \dots, n$.
We set
$\varepsilon_\alpha = \frac{1}{2}(\epsilon_{2\alpha} - \unim  
\epsilon_{2\alpha+1})$, $\alpha = 0, \dots, n-1$, and 
 $\varepsilon_{\bar \alpha} = 
\overline{\varepsilon_{\alpha}}$.
We also use the notation  $\{\epsilon^i\}$,
$\{\varepsilon^\alpha\}$ and $\{\varepsilon^{\bar\alpha}\}$
for the  dual bases of  $\{\epsilon_i\}$, $\{\varepsilon_\alpha\}$
and $\{\varepsilon_{\bar\alpha}\}$, respectively.\par
 $<,>$  is the standard Hermitian  
product of $V = \C^n$. \par
$W$ denotes the subspace
$$W = span_\C\{\varepsilon_1, \dots, \varepsilon_{n-1}\}
 = \C^{n-1}\  .$$
$(M,J)$ is always a  complex manifold with complex
structure $J$ and complex dimension $n$. 
We let  $\TM = TM \setminus \{\text{zero section}\}$ and 
$\PTM = \TM/\C^*$. \par
For any $v\in T_xM$,  $\imath_v:  T_xM \overset\sim\to\longrightarrow T_v(T_xM)$
is  the  natural isomorphism
between  $T_xM$ and $T_v(T_xM)$.
Using the maps $\imath_v$,
 any vector $w\in T_{v_o}(T_xM)\subset T(T_xM)$ 
extends to a vector field $X^{(w)}$ on  $T_xM$  by letting $X^{(w)}|_v = 
\imath_v\circ \imath^{-1}_{v_o}(w)$.
We call $X^{(w)}$ the {\it trivial extension of\/} $w$. \par
\medskip
For any $v \in T_xM$,  $J$ denotes both  
the complex structure  on  $T_xM$ and  on
$T_v(T_xM)$.
The vectors
$v^{10}$ and $v^{01}$ are the  holomorphic and anti-holomorphic
components of $v$ w.r.t. $J$
$$v^{10} = \frac{1}{2}(v - \unim  J v)\ ,\qquad
v^{01} = \frac{1}{2}(v + \unim  J v)\ .$$
The {\it dilatation
 field\/} $D$ of $T(T_xM)$  is the vector field defined as
$$D|_v = \imath_v(v)\ .$$
A {\it linear frame\/} at a point $x$
 of $M$ is an $\R$-linear 
isomorphism $u\: \R^{2n} \to T_xM$. 
A {\it complex linear frame\/} at a point $x$
 is a $\C$-linear isomorphism
$u\: \C^n\to T_xM$.
We  always 
identify a linear frame $u$  
with the corresponding basis $\{f_i\}$ in $T_xM$ where
$$f_i = u(\epsilon_i)\in T_xM\ .$$
If $u$ is complex, we denote by $u^{10}$ the corresponding
 holomorphic basis, that is 
$$u^{10} = \{e_{\alpha}
 = u(\varepsilon_\alpha) = \frac{1}{2}(f_{2\alpha} - 
\unim f_{2\alpha+1})\}\ .$$ 
For any 
linear frame $u$ on $T_xM$, the point  $x = \pi(u)$ is called {\it base point
of $u$\/}. 
\par
An {\it absolute parallelism\/} is a set of vector fields
$\{ X_1, \dots, X_{2n}\}$ which are linearly independent at all points
and, hence, which constitutes a smooth field of frames on $M$.\par
\medskip
The collection of all linear frames on $M$ is denoted by 
$L(M)$; recall that it is a $GL_n(\R)$-bundle w.r.t. the projection map $\pi$. 
The collection of all 
complex linear frames is denoted by $\LC$; it is 
a principal $GL_n(\C)$-subbundle of $L(M)$.\par
It is well-known that $\LC$ admits a unique complex structure $\hat J$  
which verifies the following two conditions: 
\roster
\item"a)" the restriction of 
$\hat J$  to the vertical subspaces of
$\LC$  coincides with the complex structure 
of $GL_n(\C)$; 
\item"b)" the projection $\pi\: \LC \to M$
is holomorphic. 
\endroster
We call $\hat J$ 
{\it the standard complex structure of $\LC$\/}.\par
\medskip
 For any subbundle $P\subset L(M)$, 
we denote by   $\theta$ its  tautological 1-form, which is 
defined  as follows.
For any   frame $u = \{f_i\}\in P\subset L(M)$ and 
 any vector $X\in T_u P$, the projected vector $\pi_*(X)$
can be written as $ \pi_*(X) =  \sum_{i=0}^{2n-1}\theta^i_u(X) f_i$
for some numbers $\theta^i_u(X)$. The tautological 1-form
$\theta$ is the $\R^{2n}$-valued 1-form 
$$\theta_u(X) = \sum_{i = 0}^{2n-1} \theta^i_u(X)\cdot\epsilon_i\ .\tag 2.1$$
If $P$ is a subbundle of $\LC$,  
 any vector $X\in T_u P$ admits the decomposition  $X = X^{10} + X^{01} = 
X^{10} + \overline{X^{10}}$,  where $X^{10}$ is the holomorphic part
of $X$ w.r.t. $\hat J$.  We denote by  $\theta^\alpha_u(X^{10})$ and 
$\theta^{\bar\alpha}_u(X^{01})$ the  components
of $\pi_*(X^{10})$ and $\pi_*(X^{01})$ w.r.t. the holomorphic
and anti-holomorphic frame $u^{10}$
and $u^{01}$, respectively. In this way two 
sets of  $\C$-valued 1-forms $\theta^\alpha$ and 
$\theta^{\bar \alpha}$ are defined at all points of $P$. 
They   are called {\it holomorphic and antiholomorphic
components of the tautological 1-form $\theta$\/}.\par
\medskip
Finally, for any $A\in  \goth{gl_n(\R)}$, we denote by $A^*$
the {\it associated fundamental
vector field\/}, that is the vector field on $L(M)$ whose flow is
$$\Phi_t^{A^*}(u) = u\circ \exp(tA)\ .$$
Since $GL_n(\C)$ acts freely and transitively on the fibers of $\LC$, 
the fundamental vector fields span any  vertical 
subspace $\Cal V_u \subset T_u\LC$. Therefore if 
$P$ is a subbundle of $\LC$ (not necessarily a principal subbundle), 
 we may  consider  the 
subspace $\g_u\subset \gl$
$$\g_u = \{\ A\in \gl\ ,\ A^*_u \in \Cal V_u\}\ .\tag2.2$$
We call $\g_u$  {\it the  algebraic 
 vertical subspace of $P$ at the point $u$\/}. Notice that
$P$ is a principal subbundle if and only if 
 $\g_u$ is a subalgebra of $\gl$  independent on 
$u\in P$. 
In this case $\g_u  =  \g = Lie(G)$, where $G$ is  the structural 
group of $P$.\par
\bigskip
\subsubhead 2.2 First properties of complex Finsler manifolds
\endsubsubhead
\medskip
\definition{Definition 2.1} A {\it complex Finsler metric\/} on 
$(M,J)$ is a continuous function
$$F\: TM \longrightarrow \R^+$$
satisfying the following properties:
\roster
\item"i)" $F$ is smooth on $\TM$;
\item"ii)" $F(u) >0$, for all $u\in \TM$;
\item"iii)" $F(\lambda u) = |\lambda|F(u)$ for all $u\in \TM$ and 
any $\lambda \in \C$.
\endroster
A {\it complex Finsler manifold \/} is  a complex manifold 
$(M,J)$ endowed with a complex Finsler metric $F$.\par
A (local) biholomorphism $f\: M\to N$ between two complex
Finsler manifolds $(M, J, F)$ and $(N, J', F')$ is called {\it
(local) biholomorphic isometry\/}
if $F'(f_*v) = F(v)$, for any $v\in TM$.
\enddefinition
Note that any complex Finsler metric $F$ is in particular
a  real
Finsler metric (see e.g. [Ca],  [Ch], [Ch1], [BC], [AP], [Sp]).\par
\medskip
For any complex Finsler manifold
$(M, J, F)$, the
{\it Finsler pseudo-sphere at a point $x$\/} is the hypersurface
$$S_x =  \{ v\in T_xM\ : F(v) = 1\} \subset T_xM\ .\tag2.3$$
\par
We say that $F$
is {\it associated with  the Hermitian metric $g$\/} if for any $v\in TM$
$$F(v) = \sqrt{g(v,v)}\ .\tag2.4$$
If this is the case, for any $x\in M$ the hermitian metric $g_x$ is recovered
from $F$ by 
$$g_x(v,w) = \frac{1}{2}\bh_{v_o}(v,w)\ ,\tag 2.5$$
where $h$ is the quadratic form defined
in the following formula (2.6) and $v_o$ is any non zero vector of $T_xM$.
Note that if $\dim_\C M = 1$, then any complex Finsler metric is associated with
a K\"ahler metric $g$. 
\medskip
{\it The quadratic form $\bh$, the cubic form $\bH$ and the 
quartic form $\gh$ of a complex 
Finsler metric $F$\/} are the following 
multilinear forms  on $T(\TM)$. Let
$X, Y, Z,W\in T_v(T_xM)$  and 
$\hat X$, $\hat Y$, $\hat Z$ and $\hat W$ be their trivial extensions. Then we set 
$$\bh_v(X,Y) = \left.\hat X\left[\hat Y\left(F^2\right)\right]\right|_v\ ;\qquad
\bH_v(X,Y,Z) = \left.\hat X\left[\hat Y\left[\hat Z\left(F^2\right)\right]
\right]\right|_v\ ;\tag 2.6$$
$$
\gh_v(X,Y,Z,W) = \left.\hat X\left[\hat Y\left[\hat Z\left[\hat 
W\left(F^2\right)\right]
\right]\right]\right|_v\ .
\tag2.7$$
Since any 
set of trivial extensions commute, it is immediate to realize that $\bh$,
$\bH$ and $\gh$ are multilinear and totally symmetric in their
arguments.\par
\medskip
In all the following, for any
 $v,w, z, t, y\in T_xM$, we will use the simplified notation
 $\bh_v(w,z)$, $\bH_v(w,z,t)$ and $\gh_v(w,z,t,y)$ in place of $\bh_v(\imath_v(w), 
\imath_v(z))$,
 $\bH_v(\imath_v(w),$ $\imath_v(z),\imath_v(t))$ and
$\gh_v(\imath_v(w),\imath_v(z),$ $\imath_v(t),\imath_v(y))$, respectively.
\medskip 
{\it The  quadratic form $\widetilde{\bh}$, the cubic form 
$\widetilde{\bH}$ and the quartic form $\widetilde{\gh}$
of a 
Finsler pseudo-sphere $S_x\subset T_xM$\/} are the 
the restrictions on $T(S_x)$ of $\bh$, $\bH$ and $\gh$.
Note that, being a  Finsler pseudo-sphere  a level set of $F^2$, it follows  
that 
$\widetilde{\bh}_v$ is equal to
$$\widetilde {\bh}_v(X, Y) = \hat X'(\hat Y'(F^2))\ ,\tag2.8$$
where  $\hat X'$, $\hat Y'$ are  
 two arbitrary vector fields, which are
tangent to $S_x$ and which coincide with  $X$ and $Y$ 
at  $v\in S_x$. A similar result holds for $\widetilde{\bH}$ and $\widetilde{\gh}$.\par
\medskip
Since  $\rho_x \= (F^2 - 1)|_{T_xM}$ is a 
defining function for the Finsler pseudo-sphere $S_x$,
we have the following immediate Lemma .\par
\bigskip
\proclaim{Lemma 2.2} Let $\D_v \subset T_v(S_x)$ be the maximal
$J$-invariant subspace 
of $T_v(S_x)$
 of a Finsler pseudo-sphere $S_x$ and let
 $\D^\C_v = \D^{10}_v + \D^{01}_v$  the corresponding
decomposition into holomorphic and anti-holomorphic subspaces. Let also
$\L_v$  the Levi form of $S_x \subset T_xM$ given by $\rho_x$
(for the definition, see e.g. (3.1) in \S 3). Then
for  any $X^{10}, Y^{10} \in  \D^{10}_v$,
$$\L_v(X^{10}, Y^{10}) = 
\widetilde{\bh}_v(X^{10}, \overline{Y^{10}})\ .\tag2.9$$
\endproclaim
\medskip
\definition{Definition 2.3}
A complex Finsler metric $F$ is called
{\it strongly pseudoconvex\/} (resp. {\it Levi non-degenerate\/})
if for any $x\in M$ the Levi form of the Finsler pseudo-sphere $S_x$
is positive definite 
(resp. non-degenerate) at all points.
\enddefinition
 Note that if  a  complex Finsler metric $F$ is strictly convex
as   real Finsler
metric   (for the definition, see f.i. 
[Ch], [AP] or [Sp]), then it is  also strongly 
pseudoconvex. The converse is obviously not  true.\par
\bigskip 
The following two Lemmata give some basic 
properties of $\bh$, $\bH$ and $\gh$.\par
\medskip
\proclaim{Lemma 2.4} Let $(M,J,F)$ be a complex Finsler manifold. 
For any $0\neq v\in T_xM$ consider the trivial extension $\hat v$ of
the vector $D_v = \imath_v(v)$ and 
let  $X$, $Y$, $Z$, $X_i\in T(T_xM)$, $i = 1, \dots, k$, some
trivially extended vector
fields. Then
\roster
\item"a)" $D\left(F^2\right)_v = 2 F^2_v$ and 
$(J D)\left(F^2\right)_v =0$;
\item"b)"
$ \left. D\left(X_1\left(X_2 \left(\dots
X_k\left(F^2\right)\dots \right)\right)\right)\right|_{v} =
(2-k) \left. X_1\left(X_2 \left(\dots
X_k(F^2)\dots \right)\right)\right|_{v}\ ;$
\item"c)"
$\left.JX_1\left(X_2 \left(\dots
X_k\left(F^2\right)\dots\right)\right)\right|_{v} 
+ \left.X_1\left(JX_2 \left(\dots
X_k\left(F^2\right)\dots\right)\right)\right|_{v} 
+ \dots  $
$$\dots + \left.X_1\left(X_2 \left(\dots
JX_k\left(F^2\right)\dots \right)\right)\right|_{v}
+ \left. JD\left(X_1\left(X_2 \left(\dots
X_k\left(F^2\right)\dots\right)\right)\right)\right|_{v}  = 0\ ;$$
\item"d)"
$D^{10}(F^2) = D^{01}(F^2) = F^2$ and  
$$\bh_v(X^{10},\hat v^{10}) = 0\ ,\quad 
\bh_v(X^{10},\hat v^{01}) = X^{10}(F^2)_v\ ;$$
\item"e)" 
$$\bH_v(X^{10}, Y^{01}, v^{10}) = \bH_v(X^{10}, Y^{01}, v^{01}) = 0\ ,\tag
2.10$$
$$
\bH_v(X^{10}, Y^{10}, v^{10}) =  - \bh_v(X^{10}, Y^{10})\ ,
\bH_v(X^{10}, Y^{10}, v^{01}) = \bh_v(X^{10}, Y^{10})\ ,
\tag 2.11$$
$$\gh_v(\hat v^{01}, X^{10}, Y^{10}, Z^{01}) = 0\ ,
\quad \gh_v(\hat v^{10},X^{10}, Y^{01}, Z^{01}) = 0\ ,\tag 2.12$$
$$\gh_v(\hat v^{10},X^{10}, Y^{10}, Z^{01}) = - 
\bH_v(X^{10}, Y^{10}, Z^{01})\ ,
\tag 2.13$$
$$ \gh_v(\hat v^{01},X^{10}, Y^{01}, Z^{01}) = - 
\bH(X^{10}, Y^{01}, Z^{01})\ .\tag 2.14$$
\endroster
\endproclaim
\demo{Proof}
Consider on  $T_xM$ the  flows 
$$\Phi_t, \Psi_t\: T_x M \to T_x M \ ,\qquad
\Phi_t(v) = e^{t}\cdot v\ ,\quad
\Psi_t(v) = e^{tJ}\cdot  v\ .\tag2.15$$
Definition 2.1 (iii) is equivalent to
$$(F^2\circ \Phi_t)(v) = e^{2t} F^2(v)\ ,\qquad
(F^2\circ \Psi_t)(v) = F^2(v)\tag2.16$$
for any $v\in \TM$ and any $t\in \R$. 
If we identify any vector  $X\in T_v(T_xM)$ with the corresponding
element  in $T_{\Phi_t(v)}(T_xM)$ and $T_{\Psi_t(v)}(T_xM)$, 
 the  differentials $\Phi_t{}_*$ and $\Psi_t{}_*$
can be written as
$$\Phi_t{}_*|_v(X) = e^t \cdot X\ ,\qquad
\Psi_t{}_*|_v(X)= e^{tJ}\cdot X\ .\tag 2.17$$
Therefore 
for any trivially extended vector fields 
 $X_i \in T(T_xM)$, $i = 1, \dots k$, 
$$e^{kt}\cdot \left[\left.X_1\left(X_2 \left(\dots
\left(X_k\left(F^2\right)\right)\dots \right)\right)\right|_{e^t v}\right] = $$
$$ =
\left.\left[\Phi_t{}_*(X_1)\left[\Phi_t{}_*(X_2)\left[\dots
\left[
\Phi_t{}_*(X_k)\left[F^2\right]\right]\dots\right]\right]\right]
\right|_{\Phi_t(v)} =$$
$$
=  e^{2t}\left.\left(X_1\left(X_2 \left(\dots
X_k\left(F^2\right)\right)\right)\right)\right|_{v}\ ,\tag2.18$$
$$\left.e^{tJ}\cdot X_1\left(e^{tJ}\cdot X_2 \left(\dots
\left(e^{tJ}\cdot X_k\left(F^2\right)\right)\dots
\right)\right) \right|_{\exp(tJ) v} = $$
$$ = 
\left. \Psi_t{}_*(X_1)\left[\Psi_t{}_*(X_2)\left[\dots
\left[
\Psi_t{}_*(X_k)\left[F^2\right]\right]\dots
\right]\right]\right|_{\Psi_t(v)}= $$
$$ =
\left.X_1\left(X_2 \left(\dots\left(
X_k\left(F^2\right)\right)\dots\right)
\right)\right|_{v}\ .\tag2.19$$
Taking the derivative at $t=0$ of (2.16), (2.18) and (2.19), 
one immediately obtains  a), b) and c). \par
d)  follows  from a), b). \par
To prove e), observe that b) and c) imply that 
$$\bH_v(V,W,v) = \hat v(V(W(F^2)))_v = D(V(W(F^2)))_v = (2-2)V(W(F^2))_v = 0\ ,$$
$$\hat v(V(W(Z(F^2))))_v = (2-3)\bH_v(V, W, Z) = -\bH_v(V, W, Z)\ ,$$
and 
$$
\bH_v(V, W, Jv) = J D(V(W(F^2)))_v = - \bh_v(JV, W) - \bh_v(V, JW)\ ,$$
$$
J D(V(W(Z(F^2))))_v = - \bH_v(JV, W, Z) - \bH_v(V, JW, Z) - 
\bH_v(V, W, JZ)\ .$$
From these identities and some straightforward computations (2.10) -  (2.14)
follow.\qed
\enddemo
\bigskip
\proclaim{Lemma 2.5} Let $(M, J, F)$ be a 
strongly pseudo-convex  Finsler manifold. Then $F$ is 
associated with an Hermitian metric $g$ if and only if one of the 
following two equivalent conditions are satisfied:
\roster
\item"i)" the cubic form $\bH$ vanishes identically;
\item"ii)" for any point $x\in M$, 
any vector $0\neq v\in T_xM$ and any two trivially extended vector
fields $X, Y \in T(T_xM)$
$$\bh_v(X^{10}, Y^{10}) = 0\ .$$
\endroster
\endproclaim
\demo{Proof} If $F$ is associated with an Hermitian metric, 
then (ii) and (2.19) are clearly satisfied. Moreover, if (ii) holds,
for any three vectors $0 \neq X, Y, Z\in T_v(T_xM)$, with trivial
extensions  $\hat X, \hat Y$ and $\hat Z$, we get 
$$\bH_v(X^{10}, Y^{10}, Z^{10}) = \hat Z^{10}(\bh(X^{10}, Y^{10}))|_v = 0\ ,$$
$$ \bH_v(X^{10}, Y^{10}, Z^{01}) =  \hat Z^{01}(\bh(X^{10}, Y^{10}))|_v = 0\ ,$$
and this implies (i). So, in order to conclude, we just need
to show that  (i) implies that
$F$ is associated with an Hermitian metric. \par
Note that if (i) holds, 
for any two trivially extended vector fields $\hat X$ and $\hat Y$,
the value of $\bh_v(\hat X, \hat Y)$ is independent on $v$. Moreover, from
(2.11), $\bh_v(\hat X^{10}, \hat Y^{10}) = \bH_v(\hat X^{10}, \hat Y^{10}, \hat v^{01}) = 0$
and hence
$$\bh_v(\hat X, \hat Y) = \bh_v(\hat X^{10}, \hat Y^{01}) + 
\bh_v(\hat X^{01}, \hat Y^{10})\ .$$
So the quadratic form $g_x$ defined by (2.5) is an Hermitian metric on $T_xM$
and, by Lemma 2.4 d), $g_x(v,v) = \bh_v(v^{10},v^{01}) = D^{10}(F^2)_v = F^2_v$.\qed
\enddemo
\bigskip
From this point on, if the opposite is not stated, by  
{\it complex Finsler metric\/} and {\it complex Finsler manifold\/}
we will mean
{\it strongly pseudoconvex complex Finsler metric\/} and {\it 
strongly pseudoconvex complex Finsler manifold\/}, 
respectively.
\bigskip
\subhead 3. The non-linear Hermitian connection of a complex Finsler manifold
\endsubhead
\medskip
\subsubhead 3.1 The H-sphere bundle of a complex Finsler manifold
\endsubsubhead
\medskip
\definition{Definition 3.1}  
An {\it Hermitianized sphere bundle 
on $(M,J)$ ({\rm or, more shortly,}  H-sphere bundle)\/} 
is a   pair $(SM, \rho)$ where: \par
\roster
\item"a)" $SM \subset TM$ is a sphere bundle;
\item"b)" each sphere 
$S_x \subset T_xM$ is a 
 strongly pseudoconvex hypersurface of $(T_xM, J)$;
\item"c)"  $\rho$ is a 
smooth  real function on $\TM$ such that 
$SM = \{ \ v \in \TM\ : \rho(v) = 0\}$.
\endroster
An H-sphere bundle $(SM, \rho)$ is called {\it circular\/} if  $SM$ is
 invariant with respect  to the linear group of transformations
$T^1 = \{e^{tJ}\ ,\ t\in \R\}$ and $\rho$ is $T^1$-invariant.\par 
Two H-spheres bundles $(SM, \rho)$,
$(SM', \rho')$ over $(M,J)$ and $(M', J')$, respectively,
are  {\it biholomorphically isometric\/}
if there exists a biholomorphism $f\: M \to M'$, such that 
 $\rho = \rho' \circ f_*|_{\TM}$
\enddefinition
\medskip
The main examples of H-sphere bundles 
are the  Finsler sphere bundles. \par
\medskip
\definition{Definition 3.2} The {\it H-sphere bundle of a 
complex Finsler manifold\/} $(M,J,F)$ 
is the pair  $(S^FM, \rho_F)$, where $S^FM$ is the  bundle 
of the Finsler spheres in $\TM$ and $\rho_F = F^2 - 1$.
\enddefinition
\medskip
Notice that $(S^FM, \rho_F)$ is always circular. Moreover
it is clear that two Finsler manifolds
are biholomorphically isometric if and
only if the corresponding H-sphere bundles are biholomorphically isometric.\par
\medskip  
\remark{Remark 3.3} Let  $S_x\subset T_xM$ be a sphere of 
a circular H-sphere bundle $(SM, \rho)$ and
 $\rho_x$ the restriction $\rho_x = \rho|_{T_xM}$, so that 
$S_x = \{\rho_x(v) = 0\}$. Let also
$\D_x\subset TS_x$ be  the family of the
 maximal  $J$-invariant tangent spaces of 
$S_x$  and 
 $\L_x$ the Levi form of $S_x$, i.e. the collection of the 
Hermitian forms on the spaces  $\D_x|_v\subset T_vS_x$, $v\in S_x$, defined by
$$\L_x(X, Y)|_{v} = X(JY(\rho_x))|_v\ .\tag3.1$$ 
By definitions,  $\D_x$ and 
$\L_x$ are   $T^1$-invariant and  each space $\D_x|_v$, $v\in S_x$,
projects isomorphically onto the  tangent space at $[v]$ of   
$$\P_x = T_xM/\C^* =  S_x/T^1\ .$$
Hence   the Levi
form  $\L_x$  induces 
an Hermitian metric $\hat \L_x$ on each tangent projective space $\P_x$. It is not difficult
to realize that $\L_x$ is indeed a K\"ahler
metric and that it  depends 
 smoothly on the point $x$ of $M$.
\endremark
\bigskip
\subsubhead 3.2 Adapted unitary 
frames of an H-sphere bundle
\endsubsubhead
\medskip
 In the next definition we introduce the concept of
 {\it adapted\/} unitary
frames  of an H-sphere bundle. In all  
formulas,  for any   
frame  $u  = \{ f_0, \dots, f_{2n-1}\}$ we use 
the  symbols
 $f_1, \dots ,
f_{2n-1}$ also to denote the vectors in 
$T_{f_0}^{10}(T_xM)$  which 
correspond to the vectors of $u$  via the natural identification map 
$\imath_{f_0}: T_xM \to T_{f_0}(T_xM)$. \par
\bigskip 
\definition{Definition 3.4}
We say that a complex linear frame $u = \{f_0, \dots, f_{2n-1}\}$  at 
$x = \pi(u)$ 
is {\it adapted to the $H$-sphere bundle $(SM, \rho)$ 
\/} if
\roster
\item"a)" $f_0 \in S_x$ and $f_1 = J f_0$;
\item"b)" the vectors $f_2, \dots, f_{2n-1}$ span 
the maximal $J$-invariant
subspace  $\D_{f_0}$ of  $T_{f_0}S_x$;
\item"c)" the holomorphic vectors  $e_1$, \dots, $e_{n-1}$ constitute
 a unitary basis for $\D_{f_0}$ with respect to  
the Levi form $\L_x$ defined by (3.1).
\endroster
\par
The subbundle  $\ULSM \subset \LC$ of all adapted unitary frames
 of $(SM, \rho)$ is 
called {\it unitary frame bundle of $(SM, \rho)$\/}. \par
If $(SM, \rho)$
is the  H-sphere bundle of 
the complex Finsler manifold $(M,J,F)$,  its unitary frame bundle 
is denoted by
$U_F(M)$.
\enddefinition
\medskip
It is immediate to realize
that $\pi\: \ULSM\to M$ is a   subbundle of $\LC$
on which $U_{n-1}$ 
 acts freely and fiber preserving. \par
Moreover, if $(SM, \rho)$ is circular, 
then for any fiber  $\Bbb V_x = \pi^{-1}(x) \subset \ULSM$, 
the quotient
$\Bbb V_x/T^1$  
is  equivalent to 
the unitary
frame bundle $U_{\hat \L_x}(\P_x)$, where $\hat \L_x$
is the K\"ahler metric defined in Remark 3.3. This implies 
$\ULSM$ is a principal bundle  only if 
all  compact K\"ahler manifolds  $(\P_x, \hat \L_x)$ are  homogeneous
spaces of  a complex subgroup $G\subset
GL_n(\C)$ which properly  contains $U_{n-1} \times T^1$. 
Since this condition is very strong, it is natural to expect 
that generically
 $\ULSM$ is  {\it not\/} a principal 
bundle.\par
In fact:\par
\medskip
\proclaim{Proposition 3.5} 
The unitary frame bundle $\ULSM\subset \LC$ of an H-sphere bundle
$(SM, \rho)$ 
is a principal subbundle  if and only if 
it is the unitary frame bundle 
of an Hermitian metric $g$ on $(M, J)$
\endproclaim
\demo{Proof} If $\ULSM$ is a principal
$G$-bundle, the group $G$ 
 verifies the following conditions: 
\roster
\item"i)" it 
is compact; 
\item"ii)" it  acts transitively 
on each sphere $S_x = \Bbb V_x/U_{n-1}$; 
\item"iii)" the isotropy subgroup of 
the $G$-action on each sphere $S_x = \Bbb V_x/U_{n-1}$ is $U_{n-1}$.
\endroster
From the 
list of the compact Lie groups 
acting transitively on a sphere ([MS], [Bo1], [Bo2]),  it follows that 
the only group which verifies i), ii) and iii) is  $G = U_n$.
By standard arguments this implies
that $\ULSM = \UG$ for some Hermitian metric. 
\qed\enddemo
\bigskip
The following Lemma gives an alternative way  to define the 
adapted unitary frames of a complex Finsler
manifold.\par
\bigskip
\proclaim{Lemma 3.6} A frame $u = \{f_i\} \in \LC$ belongs to $\UF$ if and 
only if the corresponding holomorphic frame $u^{10} = \{e_\alpha\}$
verifies 
$$\bh_{f_0}(e_\alpha, e_{\bar\alpha}) = \delta_{\alpha \beta}\tag 3.2$$
for any $0\leq \alpha, \beta \leq n-1$.
\endproclaim
\demo{Proof} By definition of adapted frame, 
$u\in \UF$ if and only if it verifies the follow three conditions  for
  $1\leq 
\lambda, \mu \leq n-1$:
$$a)\ \ F^2(f_0) = 1\ ;\quad
b)\  e_\lambda(F^2)|_{f_0} = 0\ ;
\quad c)\ \L_x|_{f_0}(e_\lambda, e_{\mu})  =
\widetilde{\bh}_{f_0}(e_\lambda, e_{\bar\mu}) = \delta_{\lambda \mu}\ .$$
Since $\left.e_0(F^2)\right|_{f_0} = \left.
D^{10}(F^2)\right|_{f_0}$, by Lemma 2.4 d) and the identity
between $\bh$ and $\widetilde{\bh}$ on the tangent spaces of the Finsler
spheres, the conditions  a), b) and c) can be rewritten in the form
$$a')\  
\bh_{f_0}(e_0, e_{\bar 0}) = 1\ ;\quad
b')\  \bh_{f_0}(e_{\lambda}, e_{\bar 0}) = 0\ ;
\quad c')\  \bh_{f_0}(e_\lambda, e_{\bar\mu}) 
=\delta_{\lambda \mu}\ ,$$
which is simply (3.2).\qed
\enddemo
\medskip
From Lemma 3.6 and Proposition 3.5, it follows that
$U_F(M)$ is a principal 
bundle if and only if $F$ is associated with an 
Hermitian metric $g$  and that in this case $\UF = \UG$.\par
\bigskip
\subsubhead 3.3 Linear and non-linear connections of Hermitian type
\endsubsubhead
\medskip
\definition{Definition 3.7} Let $P \subset \LC$ be a subbundle of $\LC$,
$\imath : P \to \LC$ be the  immersion map and $\hat J$ the standard 
complex structure of $\LC$. Then:
\roster
\item 
a {\it (non-linear) connection on 
$P\subseteq \LC$\/}
is  a   distribution $\H$ which is of real dimension $2n = \dim M$ and
complementary to the vertical distribution; 
\item if $P$ is a principal $G$-bundle, the   connection $\H$ is 
called {\it linear\/} if it is $G$-invariant; if 
$P$ is  the unitary frame bundle $\ULSM$ of a circular
H-sphere bundle $(SM,\rho)$, $\H$ is called {\it nice\/}
if it is $U_{n-1}\times T^1$-invariant;
\item 
a nice (non-linear) connection $\H$ on $\ULSM$ is called {\it isometrically 
invariant\/} if for any 
biholomorphic isometry $f$  of $(SM, \rho)$,
the lift  $\hat f$ on $\LC$ leaves $\H$ invariant;
\item 
a (non-linear) connection $\H$ is called
{\it of Hermitian type\/} if it is $\hat J$-invariant, 
i.e. for any $u\in P$
$$\hat J (\imath_*({\Cal H}_u)) = \imath_*(\Cal H_u)\ ;\tag3.3$$
\item the 
{\it connection form\/} 
 of a (non-linear) connection   $\H$ is
the unique $\gl$-valued 1-form $\omega$
on $P$, 
which vanishes on $\H$ and  verifies 
$$\omega(A^*_u) = A\tag3.4$$
for any $u\in P$ and any  $A$ in the algebraic 
vertical  subspace $\g_u \subset \gl$.
\endroster
\enddefinition
\medskip
\remark{Remark 3.8}
If $\H$ is a nice (non-linear) connection on
$\ULSM$, for any curve
$$\gamma: [a, b] \to \PTM = \ULSM/U_{n-1}\times T^1$$ 
and any frame $u\in \pi^{-1}(
\gamma_a) \subset \ULSM$, there exists a 
unique horizontal curve 
$\hat \gamma: [a,b]
\to \ULSM$, which is tangent to $\H$, projecting onto $\gamma$
and  with  $\hat \gamma_a = u$. 
Since $\H$ is  $U_{n-1}\times T^1$-invariant,
if $\hat\gamma_a^{10} = u^{10} = \{e_0, \dots , e_{n-1}\} $ and
$\hat\gamma_b^{10} = \{e'_0, \dots , e'_{n-1}\} $,  
 the linear map 
$$T_\gamma: T_xM \to T_{x'}M\ ,\qquad 
X = X^i e_i\ \  {\overset{T_\gamma}\to{\mapsto}}\ \  X^i e'_i
\tag 3.5
$$
does not dependent on the frame $u$, but only on the curve 
$\gamma$. We call it
{\it the parallel transport along $\gamma$\/}.\par
Furthermore, if $\ULSM$ is a principal subbundle of $\LC$ and $\H$ is 
a linear connection,  
the parallel transport (3.5) depends just on the curve  $\gamma_o: [a,b]\to M$
which is obtained by  projecting on $M$  the curve $\gamma$ of $\PTM$. 
In particular, 
(3.5) is the classical  parallel transport associated to a 
linear connection and it defines a linear covariant derivation on $M$.\par
If 
$\H$ is   non-linear,  the parallel transport 
(3.5)  defines the following 
 {\it non-linear covariant derivation
$\nabla$ on $M$\/}:
let $Y$ be a 
local vector field on $M$, 
$\hat X$ a vector in $T_{[v]}\PTM$  and 
$\gamma: [a,b] \to \PTM$ a curve such that $\dot \gamma_0 = \hat X$ ;  
then
$$\nabla_{\hat X} Y  \= \lim_{h\to 0} \frac{1}{h}
\left[T_{\gamma|_{[0,h]}}^{-1}
(Y_{x_h}) - Y_x\right]\ ,\tag3.6$$
where $x$ and $x_h$ are the base points of $\gamma_0 = [v]$ and $\gamma_h$,
respectively.
If we denote by 
$$\hat \pi: \ULSM \to \PTM =
\ULSM/U_{n-1}Ê\times T^1$$
 the standard projection 
map, then $\nabla_{\hat X} Y|_v$ is equal to
$$\nabla_{\hat X} Y|_v = u\left(
{\Cal X}_u\theta({\Cal  Y}) +
\omega_u({\Cal X})\cdot \theta_u({\Cal Y})\right)\ ,
\tag3.7
$$
where  $\omega$ is the connection 
form of $\H$, 
$u$ is any frame in $\hat \pi^{-1}(v)\subset 
\ULSM$  
 and ${\Cal X}$ and ${\Cal Y}$ are  two vector fields on $\ULSM$
such that $\hat \pi_*({\Cal X})_v = \hat X$ and $\pi_*({\Cal Y}) = Y$.\par
The derivation (3.7) is linear w.r.t. the vector $\hat X$, but 
{\it it is  not linearly dependent on 
the  projection of $\hat X$  in $TM$\/}, in clear contrast with 
the properties of the linear connections. 
\endremark
\bigskip
We need  now to  characterize the non-linear 
connections of Hermitian type on a given subbundle $P\subset \LC$.\par
By standard properties   of complex
fiber bundles,  $\LC$ admits at least one
linear connection  of Hermitian type. Using such connection
as a sort of  'reference point', 
we describe all   (non-linear) connections of Hermitian type 
by means of the  following proposition. \par
\medskip
\proclaim{Proposition 3.9}  Let  
$P\subset \LC$ be a subbundle of $\LC$  and denote by $\Cal F$ and
$\Cal F_P$  the families of  all non-linear 
connections of Hermitian type on $\LC$ and on $P$, respectively.\par
Let $\Cal C$ be
a linear connection  
of a Hermitian type on $\LC$ and  $\omega$ its connection form.
Then:  
\roster
\item there exists a 1-1 correspondence between
the set $\Cal G$ of smooth maps 
$E: \LC \to  \Hom_\C(\C^n, \goth{gl}_n(\C))$ and the family  $\Cal F$ 
of non-linear connections; 
for any  $E \in \Cal G$  the corresponding  
connection $\H$ is defined by 
$$
X \in \H_u\quad\Leftrightarrow\quad 
(\omega_u - E_u\circ \theta_u)(X) = 0\ ;
\tag3.8$$
\item  let $\Cal G_P$ the family  of smooth maps $E: P \subset \LC
\to \Hom_\C(\C^n, \goth{gl}_n(\C))$ such that 
$$(\omega_u - E_u\circ \theta_u)(X) = 0 \qquad
\Rightarrow\qquad X \in T_uP\ ;\tag3.9$$
then  there exists  a 1-1
map between   $\Cal G_P$ and $\Cal F_P$; for any $E\in \Cal G_P$ the 
associated connection 
is given by (3.8);
\item if $P$ is a principal $G$-subbundle (resp.  the unitary 
frame bundle of a circular H-sphere bundle),  
a  connection of Hermitian type on $P$ is a linear connection 
 (resp. a nice connection)
if and only if the corresponding map $E \in \Cal G_P$
verifies
$$E_{g\cdot u}(v,w) = g^{-1}\cdot E_u(g\cdot v,g\cdot w)\tag3.10$$
for any $g\in G$ (resp. $g\in U_{n-1} \times T^1$), $v,w\in \C^n$ and $u\in P$.
\endroster
\endproclaim
\demo{Proof}  Let $\H$ be a $\hat J$-invariant connection on 
$\LC$ and  $u$  a  frame
in $\LC$. For any $v\in \C^n$, we denote by $\hat v$ and $\hat v'$
 the  vectors
in $\Cal C_u$ and $\H_u$  so that 
$\theta(\hat v) = \theta(\hat v') = v$. It is clear that $\hat v$ and 
$\hat v'$ are uniquely determined by $u$ and $v$ and that $\hat v'
-\hat v$ is a vertical vector in $T_u\LC$. Hence there
exists a unique element $E_u(v)\in \gl$ so that
$$\hat v' -\hat v = \left.E_u(v)^*\right|_u\ .\tag3.11$$
By definitions, the map
$v\mapsto E_u(v)$ is linear for any frame $u$. 
Moreover, 
since $\Cal C$ and $\H$ are $\hat J$-invariant 
$$w = J_o v\qquad \Rightarrow\qquad\hat w = \hat J \hat v\ ,
\qquad \hat w' = \hat J \hat v'$$
and hence
$$\left.E_u(J_ov)^*\right|_u = 
\hat J\cdot \left.E_u(v)^*\right|_u = \left.(J_o\cdot E_u(v))^*\right|_u\ .$$
This shows that $E_u \in \Hom_\C(\C^n, \gl)$. Furthermore, 
for any $\hat v'\in \H_u$
$$\omega(\hat v') = \omega(\hat v + \left.E_u(v)^*\right|_u) = 
E_u(v)$$
which implies that $X\in \H_u$ if and only if 
$\omega(X) - E_u(\theta(X)) = 0$. This proves that
 any  $\hat J$-invariant 
connection $\H \in \Cal F$ determines a smooth map 
$E\in \Cal G$, which verifies  (3.8). 
It can be easily checked that  this
 map from $\Cal F$ to $\Cal G$ is 1-1.\par
(2) is proved with the same arguments used for (1).\par
(3) follows  from (3.11), the transformation law  
of the fundamental vector fields under the action of   
$GL_n(\C)$
 and the definitions of
linear and nice connections.
\qed
\enddemo
\bigskip
\subsubhead 3.4 The non-linear Hermitian connection of 
a complex Finsler manifold
\endsubsubhead
\medskip
\proclaim{Theorem 3.10} For any complex Finsler manifold $(M,J,F)$,
the unitary 
frame bundle $U_F(M)$ admits  
a unique  non-linear connection of Hermitian type.\par
This  non-linear connection of Hermitian type is nice and 
isometrically invariant.
\endproclaim
\medskip
We call the connection of Theorem 3.10 {\it the non-linear
Hermitian connection of $(M,J,F)$\/}.\par
\medskip
\demo{Proof} If $f\: M\to M$ is a local
biholomorphism and $\hat f\: L(M) \to L(M)$ is the lifted map
on the linear frame bundle, it is clear that
$\hat f$ maps $\LC$ into itself and 
is
a local biholomorphism w.r.t. to the standard complex structure
$\hat J$. 
Therefore any local biholomorphic isometry $f$  of $(M,J,F)$ is so that 
 $\hat f(U_F(M))\subset U_F(M)$ and 
 $\hat f_*$ transforms  any $\hat J$-invariant
 non-linear connection
into another $\hat J$-invariant non-linear connection. 
Hence, if there exists a unique   non-linear
connection of Hermitian type,  this connection is 
isometrically invariant. \par
\medskip
By Proposition 3.9, 
in order to prove the existence and uniqueness 
of a non-linear connection of Hermitian type it is sufficient to consider
a linear connection of Hermitian type $\Cal C$ on 
$\LC$  with connection form $\omega$, and
to show that there exists a unique map 
$E: \UF \to  \Hom_\C(\C^n, \gl)$ which belongs to  the 
family $\Cal G_{\UF}$  defined by (3.9). \par
Note that  a map $E: \UF \to \Hom_\C(\C^n, \gl)$
can be expressed in the form
$$E_u =  E^\alpha_{\beta \gamma}(u) \varepsilon_\alpha\otimes 
\varepsilon^\beta\otimes \varepsilon^\gamma + 
E^{\overline\alpha}_{\bar \beta \bar\gamma}(u) \overline{\varepsilon_\alpha}
\otimes\overline{\varepsilon^\beta}\otimes \overline
{\varepsilon^\gamma}\tag3.12$$
for some suitable $\C$-valued functions
$E^\alpha_{\beta \gamma}$,
$E^{\overline\alpha}_{\bar \beta \bar\gamma} = \overline{E^\alpha_{\beta \gamma}}
:\UF \to \C$.\par
Fix a frame $u= \{ f_i \} \in \UF$, with associated holomorphic
frame $u^{10} = \{ e_\alpha\}$, and
denote by $\hat f_i$ and 
by $\hat e_\alpha$ the unique 
vectors in $\Cal C_u$ and $\Cal C^\C_u$ such that 
$$\pi_*(\hat f_i) = f_i\ ,\qquad
\pi_*(\hat e_\alpha) = e_\alpha\ .\tag3.13$$
A map $E: \UG \to \Hom_\C(\C^n, \gl)$
verifies (3.9) at the frame $u$ if and only if for any 
$0 \leq k \leq  n-1$ the vectors
$$\hat f'_{2k} = \hat f_{2k} + (E_u(\epsilon_{2k}))^*|_u\ ,\qquad
\hat f'_{2k+1} = \hat f_{2k+1} + (E_u(\epsilon_{2k+1}))^*|_u\tag3.14$$
are tangent to $\UF$.\par
To characterize the vectors tangent to $\UF$ we can use
the next lemma, which follows immediately from 
Lemma 3.6.\par
\medskip
\proclaim{Lemma 3.11} Let $u\in \LC$, $x = \pi(u)\in M$  and $w\in T\LC$. Then 
$w \in T\UF$
 if and only if for some 
curve 
$a(t) = \{ f_i(t)\} \in \LC$, such that $a(0) = u$ and 
$\dot a(0) = w$, the associated holomorphic frames $a^{10}(t) = \{e_\alpha(t)\}$
and the curve of vectors $f_0(t)$ verify
$$\bh_{f_0(0)}(e_\alpha(0), 
e_{\bar \beta}(0)) = \delta_{\alpha\beta}\ ,\qquad
\left.
\frac{d}{dt}\left[\bh_{f_0(t)}(e_\alpha(t), 
e_{\bar \beta}(t))\right]\right|_{t=0} = 0\ .
\tag3.15$$
\endproclaim
\medskip
For any $i = 0, \dots, 2n-1$, let us 
consider a curve $g_{i}\: (-1, 1) \subset \R\to \LC$ and a curve
$d_{i}\: (-1, 1) \subset \R\to GL_n(\C)$ so that 
$$a_{i}(0) = u\ , \quad
\dot a_{i}(0) = \hat  f_i\ ,\qquad 
g_{i}(0)= \Bbb I\ ,\quad \dot g_{i}(0) = 
E_u(\epsilon_{i})\ .\tag3.16$$
Note that, by definitions,  the curve $a'_{i}(t) = 
g_{i}(t)\cdot a_{i}(t)$ verifies 
$$a'_{i}(0) = u\ ,\quad
\dot a'_{i}(0) =  \hat f'_i\ ,\quad \pi(a'_{i}(t)) =
\pi(a_{i}(t))\ .\tag
3.17$$
For any $i = 0, \dots, 2n-1$, we will denote  
by $f^{(i)}_0(t)$ and by $f^{(i)}_0{}'(t)$ 
the curve of vectors determined by the first 
vectors of the frames $a_i(t)$ and $a'_i(t)$, respectively;
by $e_\alpha^{(i)}(t)$ and by $e_\alpha^{(i)}{}'(t)$ we denote the 
vectors of the holomorphic frames  $a^{10}_i(t)$
and $a^{10}_i{}'(t)$. We will also use the notation
$$\bh_{\alpha  \bar\beta, i}(t) = \bh_{f^{(i)}_0(t)}(e^{(i)}_\alpha(t),  
e^{(i)}_{\bar \beta}(t))\ , \qquad
\bh'_{\alpha  \bar\beta, i}(t) = \bh_{f^{(i)}_0{}'(t)}(e^{(i)}_\alpha{}'(t),  
e^{(i)}_{\bar \beta}{}'(t))\ ,$$
$$H_{\alpha \bar \beta \zeta, i}(t) = 
\bH_{f^{(i)}_0(t)}(e^{(i)}_\alpha(t),  
e^{(i)}_{\bar \beta}(t), e^{(i)}_{\zeta}(t))\ ,$$
$$H_{\alpha \bar \beta \bar\zeta, i}(t) = 
\bH_{f^{(i)}_0(t)}(e^{(i)}_\alpha(t),  
e^{(i)}_{\bar \beta}(t), e^{(i)}_{\bar\zeta}(t))\ ,$$
$$H_{\alpha \bar \beta \zeta} = 
\bH_{f_0}(e_\alpha,  
e_{\bar \beta}, e_{\zeta})\ ,\qquad
H_{\alpha \bar \beta \bar\zeta} = 
\bH_{f_0}(e_\alpha,  
e_{\bar \beta}, e_{\bar\zeta})\ .$$
Note that for $i = 2\gamma$ or $i = 2\gamma + 1$, 
$\gamma = 0, \dots, n-1$, we have
$$\left.\frac{d}{dt}\left(\bh'_{\alpha \bar\beta, 2\gamma}\right)\right|_{t=0} = 
\left.\frac{d}{dt}\left(\bh_{\alpha  \bar\beta, 2\gamma}\right)\right|_{t=0}
+ E^{\delta}_{\alpha\gamma}\bh_{\delta  \bar\beta,2\gamma}(0)
+  E^{\bar\delta}_{\bar\beta\bar\gamma}
\bh_{\alpha \bar\delta,2\gamma}(0) + $$
$$+
E^{\zeta}_{0 \gamma}
H_{\alpha \bar \beta \zeta, 2\gamma}(0)
+ E^{\bar\zeta}_{\bar 0 \bar\gamma}
H_{\alpha \bar \beta \bar\zeta, 2\gamma}(0)
= $$
$$  = \left.\frac{d}{dt}\left(\bh_{\alpha \bar\beta, 2\gamma}\right)\right|_{t=0}
+ E^{\beta}_{\alpha\gamma}
+  E^{\bar\alpha}_{\bar\beta\bar\gamma}
+ 
E^{\zeta}_{0 \gamma}
H_{\alpha \bar \beta \zeta}
+ E^{\bar\zeta}_{\bar 0 \bar\gamma}
H_{\alpha \bar \beta \bar\zeta}\tag3.18$$
$$\left.\frac{d}{dt}\left(\bh'_{\alpha  \bar\beta, 2\gamma
+ 1}\right)\right|_{t=0} = 
\left.\frac{d}{dt}\left(\bh_{\alpha  \bar\beta, 2\gamma+1}\right)\right|_{t=0}
+\unim  E^{\beta}_{\alpha\gamma}
- \unim  E^{\bar\alpha}_{\bar\beta\bar\gamma} + $$
$$ + \unim 
E^{\zeta}_{0 \gamma}
H_{\alpha \bar \beta \zeta}
- \unim   E^{\bar\zeta}_{\bar 0 \bar\gamma}
H_{\alpha \bar \beta \bar\zeta}\tag3.19$$
where we used the fact that, for any $i = 0,\dots, 2n-1$, $a^{(i)}(0) = u$
and hence that
 $f^{(i)}_0(t) = f_0$  and  that
 $\bh_{\delta  \bar\beta,i}(0) = \bh_{f^{(i)}_0(0)}
(e_\delta(0),e_{\bar \beta}(0)) = 
\delta_{\delta \beta} =  \bh'_{\delta \bar\beta,i}(0)$.\par
From Lemma 3.11,  $\hat  f_i'$ is tangent to $\UF$ if and only if
$\left.\frac{d}{dt}\left(\bh'_{\alpha \bar \beta,i}(t)\right)\right|_{t=0} = 0$
for any $i = 0, \dots, 2n-1$. \par
Note that, by (2.10), $H_{0  \bar \beta \zeta} = 0$. Therefore, 
if 
$\alpha = 0$,
 (3.18) and (3.19) vanish if and only if 
$$E^{\beta}_{0\gamma} = 
- \left.\frac{d}{dt}\left(\bh_{0 \bar \beta, 2\gamma}\right)\right|_{t=0}
+\unim  
\left.\frac{d}{dt}\left(\bh_{0 \bar \beta, 2\gamma+1}\right)\right|_{t=0}
\ .\tag
3.20$$
By substituting (3.20)  back into (3.18) and (3.19), we obtain
 all others components $E^\beta_{\alpha\gamma}$. In fact,  
$$E^{\beta}_{\alpha\gamma} = 
- \left.\frac{d}{dt}\left(\bh_{\alpha \bar \beta, 2\gamma}\right)\right|_{t=0}
+\unim  
\left.\frac{d}{dt}\left(\bh_{\alpha \bar\beta, 2\gamma+1}\right)\right|_{t=0}
- 2 E^{\zeta}_{0 \gamma}
H_{\alpha \bar \beta \zeta}
\tag3.21$$
This implies that  $\Cal G_{\UF}$ contains exactly one map $E$.  
 From (3.20) and (3.21), it can  be 
checked that $E$ verifies (3.10) for any element 
$g\in T^1 \times U_{n-1}$
and hence that it defines a nice non-linear connection by Proposition 3.9 (3).
\qed
\enddemo
\bigskip
Since the unitary frame bundle 
$\UG$ of an 
Hermitian metric $g$ on $(M,J)$ coincides with the unitary 
frame bundle $U_F(M)$ of the Finsler metric $F(v) = \sqrt{g(v,v)}$, 
 Theorem 3.10 gives the following
classical result  as an immediately corollary (see e.g. [KN], vol. II): 
{\it an Hermitian manifold $(M, J,g)$ has  a unique  Hermitian
linear connection\/}.\par
\bigskip
\subhead 4. The absolute parallelism on $\UF$ and Kobayashi's theorem
\endsubhead
\medskip
\subsubhead 4.1 Pseudo-hermitian structures on a real hypersurface
\endsubsubhead
\medskip
In this subsection, we recall the definition of
pseudo-hermitian structure on CR manifold of codimension one and 
Webster's theorem on the existence and uniqueness of an invariant 
linear connection for any pseudo-hermitian structures. This result 
is essential  for the construction of an invariant absolute parallelism
on the unitary frame bundle $\UF$ of a complex Finsler manifold.\par
\medskip
Let $S$ be a $(2n-1)$-dimensional manifold. 
{\it A
 CR structure on $S$\/} is a pair $(\D, J)$, where 
$\D\subset TS$ 
is a  distribution and $J$ is a 
smooth family of complex 
structures $J_p$ on the subspaces $\D_p \subset T_pM$. It is called
{\it integrable\/} if the
 holomorphic  distribution  $\D^{10} \subset \D^\C$ defined by $J$
is closed under Lie brackets.  It  is called
of {\it codimension p\/} if the distribution $\D$ is of codimension $p$.\par
An  CR structure $(\D, J)$ of codimension one 
is called {\it  Levi non-degenerate\/}
if $\D$  is 
a contact distribution,  i.e.  if for any local
1-form $\theta$ such that $\ker \theta = \D$, then
 $d\theta_p$ 
is non degenerate on $\D_p$ at any point $p$ where $\theta$ is defined.
\medskip
\definition{Definition 4.1}[We] A {\it pseudo-hermitian structure\/} on $S$
is a pair $((S, \D, J); \theta)$ where $(S, \D, J)$ is a codimension
one Levi non-degenerate
 CR structure and $\theta$ is a 1-form
 on  $S$  such that $\ker \theta_p = \D_p$ for any $p\in S$.
\par
A {\it pseudo-hermitian transformation\/} of $((S, \D, J); \theta)$
is a diffeomorphism $f: S\to S$ such that  $f_*(\D) \subset \D$,
$f_* \circ J|_\D = J\circ f_*|_\D$ and $f^*\theta = \theta$.
\enddefinition
In the following,  we denote a pseudo-hermitian structure 
only by a pair $(S, \theta)$.\par
\medskip
A standard example of pseudo-hermitian structure is the 
following. Let $S$ be a smooth real hypersurface
in $\C^n$ and  $(\D, J)$ the codimension one CR structure determined by 
the maximal $J_o$-invariant subspaces in $TS$ and the complex
structures determined by the complex structure  $J_o$ of $\C^n$. Assume 
also that $\rho$ is a smooth  defining function for $S$, i.e. 
$S = \{ p\in \C^n\ : \ \rho(p) = 0\}$.
The 1-form $\theta^\rho$ 
$$\theta^\rho_p(v) = d\rho_p(J_o v)\tag4.1$$
vanishes exactly on the vectors on $\D$. Hence if $(S, \D, J)$ 
is Levi non-degenerate, then 
$(S,  \theta^\rho)$ 
is a pseudo-hermitian structure.
\medskip
\definition{Definition 4.2}
Let $(S,\theta)$ be a pseudo-hermitian
structure and  let $u = \{f_1, $ $\dots, f_{2n-1}\}$  
a linear frame at a point   $p\in S$.
The frame $u$ is called {\it adapted to $(S, \theta)$ \/} if 
\roster
\item"a)" $\theta(f_1) = 1$ and $\theta(f_i) = 0$ for 
$2\leq i\leq 2n-1$;
\item"b)" $d\theta(f_i,f_j) = \delta_{ij}$ for $2\leq i,j\leq 2n-1$;
\item"c)" $J f_{2i} = f_{2i + 1}$ for $1 \leq i \leq n-1$.
\endroster
The collection  $U_\theta(S)$ of all adapted frames of frames 
of a pseudo-hermitian
structure $(S, \theta)$ is called {\it  unitary frame bundle of 
$(S, \theta)$\/}.
\enddefinition
\medskip
Conditions b) and c) can be restated claiming that the vectors
$e_\alpha = f_{2\alpha} - i f_{2\alpha +1}$, with $\alpha = 1, \dots, n-1$,
constitute a  holomorphic basis for $\D^{10}_p \subset \D^\C_p$, which is
unitary  w.r.t. the Levi form $\L(X, Y) \= d\theta_p(X, \bar Y)$. 
It can  be checked that 
$U_\theta(S)$ is 
a principal subbundle of  $L(S)$ with 
structural group $U_{q,q'}$, where $(q, q')$ is the signature
of the Levi form $\L$ (see also [We]).\par
\medskip
For any linear connection $\H$  on $U_\theta(S)$ and  any frame $u
 = \{ f_1, \dots, f_{2n-1}\} 
\in U_\theta(S)$,
let us denote by
$p_u = \pi_*|_{\H_u}: \H_u \to T_{\pi(u)}S$ the 
restriction  of 
$\pi_*$ on $\H$. We also denote by  $\{ \tilde f_i = p^{-1}_u(f_i)\}$
the  basis of $\H_u$ which projects onto the vectors of $u$. \par
Since $p_u$ is 
a linear isomorphism between $\H_u$ and
$T_{\pi(u)}S$,  we may always consider the subspace $\tilde \D_u = 
p_u^{-1}(\D_{\pi(u)}) \subset \H_u$
and the complex structure $\tilde J$ on  $\tilde{\D}_u$, defined by
$$\tilde J\cdot v = (p_u^{-1}\circ J \circ p_u)(v)\ .\tag4.2$$
$(\tilde \D, \tilde J)$ are called {\it the
horizontal lifts of the CR structure $(\D,J)$\/}.
\medskip
\definition{Definition 4.3} Let $\H$ be a linear connection on 
$U_\theta(S)$ and let $(\tilde \D, \tilde J)$ the corresponding horizontal
lift of the CR structure of $S$. Let also 
$\tilde e_\alpha$  the holomorphic vector fields in $\tilde \D^{10}$ defined by 
$$\tilde
e_\alpha|_u = 
\tilde f_{2\alpha}|_u - i \tilde  f_{2\alpha +1}|_u\in \tilde{\D}^{10}_u\ ,\quad
1\leq \alpha\leq n-1\ .
$$
We say that $\H$ is of {\it Webster type\/} if:
\roster
\item"a)" any  Lie bracket
between two vector fields of the  holomorphic distribution
$\tilde{\D}^{10}\subset \tilde \D^\C$ is trivial;
\item"b)"  $\pi_*([\tilde e_\alpha, \overline{\tilde
e_\beta}]) = -   \unim   \delta_{\alpha\beta} \pi_*(\tilde f_1)$
for any $1\leq \alpha, \beta \leq n-1$;
\item"c)" for any $1 \leq \alpha \leq  n-1$
the complex vector field $T_\alpha = [\tilde e_{\alpha},\tilde f_1]$
is so that $\pi_*(T_\alpha)$ takes  values in 
$\D^{01} = \overline{\D^{10}}$ at all points.
\endroster
We call the vector field $T_\alpha = [\tilde e_{\alpha},\tilde f_1]$ 
{\it the $\alpha$-th component of the torsion of $\H$\/}.
\enddefinition
\medskip
We conclude with the following important result by S. Webster.\par
\bigskip
\proclaim{Theorem 4.4} [We] If 
$(S, \theta)$ is a  pseudo-hermitian structure, there
exists exactly one linear connection
of Webster type  on $U_\theta(S)$. This connection 
 is invariant under the group $Aut(S,\theta)$ of pseudo-hermitian
transformations of $(S,\theta)$.\par
Moreover, the $\alpha$-th components of the torsion
 $T_\alpha = [\tilde e_{\alpha},\tilde f_1]$ of the connection of 
Webster type 
 vanish identically   if and only if 
the vector field $f_1$, given by 
the first vectors of all adapted frames,  is
an infinitesimal CR transformations of $S$. 
\endproclaim
We call such unique linear connection  {\it the Webster connection 
of $(S, \theta)$\/}.\par
\bigskip
\subsubhead 4.2 The generalized fundamental vector
fields on the unitary frame bundle of an H-sphere bundle
\endsubsubhead
\medskip
Consider a circular 
H-sphere bundle $(SM, \rho)$. For any $x\in M$, 
let  $\rho_x$  the restriction $\rho_x = 
\rho|_{T_xM}$ and  $\theta_x = \theta^\rho_x$ the 1-form on 
 $S_x$ defined by (4.1). By definitions, 
each pair $(S_x, \theta_x)$ is a pseudo-hermitian
structure. \par
If $u = \{f_0, f_1, \dots, f_{2n}\} \subset T_xM$
 is an adapted unitary frame of $\ULSM$ and if we identify the vectors
$\{f_1, \dots, f_{2n}\}$ with the corresponding vectors at 
$T_{f_0}S_x$, we immediately see that $u$ is also an adapted frame for  
the pseudo-hermitian structure $(S_x, \theta_x)$. In other words, 
the fiber
$\Bbb V_x = \pi^{-1}(x) \in \ULSM$ can be identified 
with the unitary frame bundle $U_{\theta_x}(S_x)$.\par
 Let 
$\W_x$ the Webster connection on $\Bbb V_x \simeq U_{\theta_x}(S_x)$.
Notice that the flow of the vector field $f_1$
on $S_x$ coincides with 1-parameter
group of transformations  given by $T^1$. Therefore by Theorem 4.4
each component of the torsion of the Webster connection 
$\W_x$ vanishes identically.
\par
\medskip
Using the Webster connection $\W_x$
 and the fundamental vector fields associated
to the   Lie algebra $\u = Lie(U_{n-1} \times T^1)$ we 
define an absolute parallelism on any fiber $\Bbb V_x \subset \ULSM$
as follows.\par
\medskip
\definition{Definition 4.5} Let $\ULSM$ be the  unitary frame bundle
of a circular H-sphere bundle $(SM, \rho)$. For any $x\in M$, let also
 $\W_x$ the Webster connection on the fiber $\Bbb V_x = 
\pi^{-1}(x) \simeq U_{\theta_x}(S_x) = U_\theta(S_x)$
and  $\hat \pi: 
\Bbb V_x  \to S_x = \Bbb V_x/U_{n-1}$  the
standard projection map.
\par 
For any element
$X \in W \oplus(\u)$ we associate a  vertical 
vector field $\tilde X$ of $\ULSM$ associated  as follows:\par
\roster 
\item if $X \in \u$, we set $\tilde X = X^*$;
\item if $X \in W (= \C^{n-1})$, we set $\tilde X$
as the vector field so that, for any frame $u$, 
 $\tilde X_u$ is the unique vector in
 $(\W_{\pi(u)})_u$, such that 
$\hat \pi_*(\tilde X_u) = u(X) \in T_{f_0} S_x$;
\item for any two vectors $X, X' \in W \oplus(\u)$, we set 
 $\widetilde{X + X'} = \tilde X + \tilde X'$.
\endroster
Any vertical vector field $\tilde X$ is called {\it  
generalized fundamental vector field\/}.
\enddefinition
\medskip
If $\{\epsilon_i,  E_j\}$ is a basis for $W \oplus(\u)$, 
the generalized fundamental vector fields $\{\tilde \epsilon_i,\tilde E_i\}$
are linearly independent 
at all points and  they give an absolute 
parallelism on each fiber  $\Bbb V_x$. By construction, the 
vector fields $\{\tilde \epsilon_i,\tilde E_i\}$ are mapped into
themselves  by any biholomorphic isometry of $(SM, \rho)$.\par
One may also check that 
if $\ULSM$ is a principal subbundle (and hence with structural group $U_n$;
see Proposition 3.5), 
the  generalized fundamental vector fields coincide with the fundamental
vector fields of the elements in $\goth u_n$. \par
\bigskip
\subsubhead 4.3 The absolute parallelism on $\UF$
 and the isometry group of $(M,J,F)$
\endsubsubhead
\medskip
Let $\UF$ be the unitary frame bundle of a complex Finsler manifold
$(M,J, F)$ and  $\H$ the  non-linear Hermitian connection 
on $\UF$. \par
At any frame
$u = \{f_0, \dots, f_{2n-1}\}\in \UF$ we denote by $\hat f_i|_u$
the unique vectors in $\H_u$ which project onto $f_i \in T_{\pi(u)} M$. 
Let also $\{\epsilon_i\}$ the standard 
basis of $W = \C^{n-1}$, $t$ be a generator of $\R = Lie(T^1)$ and 
$\{E_j\}$  a basis for $\goth u_{n-1}$. Let also  $\tilde \epsilon_i$, 
$\tilde t$  and 
$\tilde E_j$ 
the corresponding generalized fundamental vector fields on $\UF$.\par
Then the set of 
vector fields  
$$\sigma^\H = \{\hat f_i, \tilde \epsilon_j, \tilde t, 
\tilde E_k\}\tag4.3$$
is an absolute parallelism
on $\UF$. It is unique, up to a change of the generator
$t$ and of the basis $\{E_j\}$  for 
$\goth u_{n-1}$ and it is invariant under all  
 local biholomorphic isometries  of $(M,J,F)$. \par
We call 
$\sigma^\H$  
{\it the absolute parallelism associated to the Hermitian connection $\H$\/}.\par
As a consequence of 
Kobayashi's theorem on the automorphism group of an absolute parallelism
([Ko]), the following result is easily obtained.\par
\medskip
\proclaim{Proposition 4.6} Let $(M, J,F )$ be a complex Finsler 
manifold of complex dimension $n$ and let $\H$ the non-linear 
Hermitian connection of $(M, J, F)$.
\roster
\item The local
holomorphic isometries of $(M, J, F)$ are in 1-1 correspondence
with the local diffeomorphisms of $\UF$ which preserve the 
absolute parallelism $\sigma^{\H}$. For any  local biholomorphic
isometry  $f$  the corresponding local diffeomorphism 
 is the restriction on $\UF\subset \LC$  of the   diffeomorphism
$$\hat f: L(M) \to L(M)\ , \qquad \hat f(u) = f_* \circ u\ .$$
\item The group $Iso_F(M)$ of all biholomorphic isometries
of $(M,J,F)$ is a Lie group of dimension less or equal to
$$\dim_\R \UF = \dim_\R V + \dim_\R W + \dim_\R (\u) = n^2 + 2n\ .$$
\item $\dim_\R Iso_F(M) = n^2 + 2n$ if and only if $F$ is 
 associated with a K\"ahler metric $g$ 
and $(M,J,g)$ is
$(\CP^n, g_{c^2})$, $(\C^n, g_0)$ or $(B^n , g_{-c^2})$, where 
$g_{c^2}$, $g_0$ and $g_{-c^2}$ denote 
the K\"ahler metrics with constant holomorphic 
sectional curvature $c^2$, $0$ and $-c^2$, respectively.
\endroster
\endproclaim
\demo{Proof} (1) It is proved  with the same 
arguments  of  Proposition 3.3 in [Sp]; they are very similar
 to the arguments of the proof of Proposition VI. 3.1 in 
[KN] Vol. I.\par
(2) It is an immediate corollary of (1) and of Theorem I.3.2 in [Ko1].\par
(3) If $\dim_\R Iso_F(M) = n^2 + 2n$, 
then $U_F(M)$ is a principal subbundle of $\LC$ and 
 any isotropy subgroup $Iso_F(M)_x$ acts
transitively on any fiber $\Bbb V_x$ of $\UF$.  By Proposition 3.5 and Lemma 3.6
it follows that $Iso_F(M)_x \simeq  U_n$ for any  $x\in M$ and 
 that  $F$ is associated with an Hermitian metric $g$. Since 
for any  $x\in M$ we have that $Iso_F(M)_x = Iso_g(M)_x \simeq  U_n$ , 
we also have that $(M,J,g)$ is an
Hermitian symmetric space,  $g$ is 
K\"ahler and the holomorphic sectional curvature is constant.\par
We claim now that $M$ is  simply 
connected. Suppose not and let $\pi: \tilde M\to M$  its universal covering map 
and $\tilde g = \pi^*g$. For any $x\in M$ and any $y\in \pi^{-1}(x)$, 
the isotropy subgroup  $Iso_{g}(M)_{x}\simeq U_n$ is embedded into 
$Iso_{\tilde g}(\tilde M)_{y}$ and hence 
  $Iso_{\tilde g}(M)_{y} = Iso_{g}(M)_{x}\simeq U_n$. This implies that 
any deck transformation $\Gamma$ belongs to the normalizer  
$N_{Iso_{\tilde g}(\tilde M)}(Iso_{\tilde g}(\tilde M)_y) $ of the subgroup
$Iso_{\tilde g}(M)_{y}$ in $Iso_{\tilde g}(\tilde M)$. 
Suppose now that $\tilde g$ has positive holomorphic sectional 
curvature. By
 the classification of simply connected complex space forms, $\tilde M$ 
is $\C P^n$, $Iso_{\tilde g}(\tilde M) = SU_{n+1}$ and 
$N_{Iso_{\tilde g}(\tilde M)}(Iso_{\tilde g}(\tilde M)_y) = 
N_{SU_{n+1}}(U_n) = U_n$. This means that any deck transformation $\Gamma$
belongs to $Iso_{\tilde g}(M)_{y}$, and this cannot be because $\Gamma$ fixes no
point. 
Suppose then that $\tilde g$ has non-positive   holomorphic sectional 
curvature. In this case for any non-trivial deck transformation and 
any $h\in Iso_{\tilde g}(M)_{y}$, there exists some 
 $h' \in 
Iso_{\tilde g}(M)_{y}$ so that $h\circ \Gamma = \Gamma \circ h'$ and hence
 $h(\Gamma(y)) = \Gamma(h'(y)) = \Gamma(y)$. 
Since $h$ fixes 
$y$ and $\Gamma(y)$, it fixes point by point 
 the unique length minimizing geodesic between 
$y$ and $\Gamma(y)$. But this cannot be 
because $Iso_{\tilde g}(\tilde M)_y = U_n$ and it
fixes no vector in $T_y\tilde M$. \par
Being $M$ simply connected, the claim follows  from 
the classification of simply connected complex space forms. \qed
\enddemo
\bigskip
\subhead 5. The  invariants
 of a complex
Finsler manifolds
\endsubhead
\medskip
\subsubhead 5.1 Notation
\endsubsubhead
\medskip
In all the following sections,
the greek indices $\alpha$, $\beta$, $\gamma$, $\delta$, $\varepsilon$  always 
run between $0, \dots, n-1$; the indices
$\lambda, \mu, \nu, \rho, \sigma$   run between $1$ and $n-1$.\par
We  denote by 
$E^\alpha_\beta = \varepsilon_\beta\otimes \varepsilon^\alpha$ the
elements of 
the standard basis of $\gl$.
An element  $A = 
A^\alpha_\beta E^\beta_\alpha\in \gl$ can be also  
expressed using just the complex matrix $A^\alpha_\beta$. \par
 For any adapted unitary
frame $u = \{f_i\}$ and corresponding holomorphic frame
$u^{10} = \{e_\alpha\}$, we set
$$h_{\alpha \beta}(u) = \bh_{f_0}(e_\alpha, e_\beta)\ ,\ 
H_{\alpha \beta \gamma}(u) = \bH_{f_0}(e_\alpha, e_\beta, e_\gamma)\ ,\  
H_{\alpha \beta \gamma\delta}(u) = \gh_{f_0}(e_\alpha, e_\beta, e_\gamma, e_\delta)
.$$
Analogous meanings will have the symbols 
$h_{\bar \alpha \bar \beta}(u)$, $H_{\alpha \beta \bar\gamma}(u)$,
$H_{\alpha  \bar \beta \bar \gamma}(u)$, etc.. 
\par
\medskip
On $\UF$ we have  the 
following distributions and CR structures:
\roster
\item"-" $\H$ is the non-linear Hermitian 
connection;
\item"-"  $\W = \bigcup_{x\in M} \W_x$ is the  
distribution  obtained  as  union of the 
 Webster connections
$\W_x$ of the fibers $\Bbb V_x = \pi^{-1}(x) = U_{\theta_x}(S_x)$;  
\item"-" $(\tilde \D, \tilde J)$ is the CR structure given by  
$\tilde \D = \bigcup_{x\in M}
\tilde D_x$ and  $\tilde J = \bigcup_{x\in M} \tilde J_x$, where each
$(\tilde \D_x, \tilde J_x)$ is the horizontal lift in $\W_x$ of the 
CR structure of the Finsler sphere $S_x$ (see Definition 4.3).
\endroster
Notice that
$(\H, \hat J)$ and 
 $(\tilde \D, \tilde J)$ are  both {\it integrable\/} CR structures on $\UF$
(see \S 4.1).\par 
\smallskip
$\omega$ is the connection form of $\H$ (see Def. 3.7) and we define
$\omega^\alpha_\beta$
as the $\C$-valued 1-forms on $\UF$ which verify 
$$\omega =  \sum_{\alpha ,\beta}
E^\alpha_\beta\otimes \omega^\beta_\alpha \ .$$
We also set $\omega^{\bar \alpha}_{\bar\beta} = \overline{\omega^\alpha_\beta}$.
We call $\omega^\alpha_\beta$ and $\omega^{\bar \alpha}_{\bar\beta}$
{\it the holomorphic and anti-holomorphic components of the
connection 
form $\omega$\/}.\par
\smallskip
We denote by $\E^\alpha_\beta$ and by $\E^{\bar \alpha}_{\bar \beta}$
the complex vector fields on $\LC$ defined by 
$$\E^\alpha_\beta = \frac{1}{2}\left[
(E^\alpha_\beta)^* - \unim  (J_o E^\alpha_\beta)^*\right]\ ,
\qquad \E^{\bar \alpha}_{\bar \beta} = \overline
{\E^\alpha_\beta}\ .$$
Note that, if we extend $\C$-linearly the 1-forms $\omega^\alpha_\beta$ and
 $\omega^{\bar \alpha}_{\bar\beta}$,
we have that
$$\omega^\alpha_\beta(\E^\varepsilon_\gamma) = \delta^\alpha_{\gamma} 
\delta^\varepsilon_\beta
\ , \qquad \omega^{\bar \alpha}_{\bar \beta} (\E^\delta_\gamma) = 0\ .$$
Moreover, we recall that by the properties  of fundamental vector fields (see [KN])
$$\Cal L_{(E^\alpha_\beta)^*} \theta = - \sum_\gamma
 [E^\alpha_\beta,
\varepsilon_\gamma]\otimes\theta^\gamma - \sum_\gamma[E^\alpha_\beta,
\varepsilon_{\bar\gamma}]\otimes\theta^{\bar\gamma}\ ,$$
$$
\Cal L_{(J_oE^\alpha_\beta)^*} \theta = -\sum_\gamma
 \unim  [E^\alpha_\beta,
\varepsilon_\gamma]\otimes \theta^\gamma  + \unim  \sum_\gamma[E^\alpha_\beta,
\varepsilon_{\bar\gamma}]\otimes\theta^{\bar\gamma}\ .$$
This implies that
$$\Cal L_{\E^{\alpha}_\beta}\theta = \frac{1}{2}
\left[\Cal L_{(E^\alpha_\beta)^*} \theta - \unim  
\Cal L_{(J_oE^\alpha_\beta)^*} \theta\right] = 
 - \varepsilon_\beta\otimes \theta^\alpha\ ,
\qquad \Cal L_{\E^{\bar\alpha}_{\bar\beta}}\theta = 
\overline{\Cal L_{\E^{\alpha}_\beta}\theta} = 
  - \varepsilon_{\bar\beta}\otimes \theta^{\bar\alpha}\ .$$
\medskip
For any fiber $\Bbb V_x = \pi^{-1}(x) \subset \Bbb V_x$, 
we denote by 
$$\hat \pi_x: \Bbb V_x = U_{\theta_x}(S_x) \to S_x = \Bbb V_x/U_{n-1}$$
the standard projection map of  $\Bbb V_x = U_{\theta_x}(S_x)$
onto $S_x =  \Bbb V_x/U_{n-1}$. \par
For any $i = 0, \dots, 2n-1$,  $\hat f_i$ is the 
vector field in $\H$ such that at all $u\in \UF$
$$\pi_*(\hat f_i)|_u = f_i = u(\epsilon_i) \in T_{\pi(u)} M\ .$$
For any $a= 2, \dots, 2n-1$, 
$\tilde \epsilon_a$ is
 the generalized fundamental vector fields which corresponds
to elements of the real basis $\{ \epsilon_2, \dots, \epsilon_{2n-1}\}$
of $W = \C^{n-1}$ (see \S 2.1). By construction,
any vector field $\tilde \epsilon_a$ is vertical, it  takes values in
$\tilde \D\subset \W$ and $\tilde J \tilde \epsilon_{2i} = \tilde
\epsilon_{2i+1}$.\par
\smallskip
Finally, we denote by $\hat e_\alpha$ and $\tilde e_\lambda$ the holomorphic 
vector fields
$$\hat e_\alpha = \frac{1}{2}\left(\hat f_{2\alpha} - \unim  
 \hat f_{2\alpha + 1}\right)\ ,\qquad
\tilde e_\lambda = \frac{1}{2}\left(\tilde \epsilon_{2\lambda} - 
\unim  \tilde \epsilon_{2 \lambda +1}\right)
\ .$$
They  coincides with the complex 
vector fields in $\H^\C$ and $\tilde \D^\C \subset \W^\C$, which 
 are  mapped by $\pi_*$ and $\hat \pi_x{}_*$ 
onto the holomorphic vectors of the corresponding adapted frames of $T_xM$
and of $T_{\hat \pi_x(u)}S$, respectively.
\bigskip
\subsubhead 5.2 The algebraic vertical subspaces of $\UF$ and 
the distribution $\W$
\endsubsubhead
\medskip
We want to determine the algebraic vertical subspaces
of $\g_u$ (see definition in \S 2.1). For this
we give the following technical lemma, which follows directly from
definitions and Lemma 2.4 e).\par
\medskip
\proclaim{Lemma 5.1}  For any choice of the indices $A, B, C$, the 
functions $H_{A B C}$ are totally symmetric w.r.t. to $A$, $B$ and $C$. 
Furthermore, for any $u = \{ f_0, \dots, f_{2n-1}\}$, 
$$H_{\alpha \bar \beta 0}(u)  = H_{\bar \alpha \beta \bar 0}(u) 
= 0\ , \ \  
H_{\alpha \beta 0}(u) = - h_{\alpha \beta}(u)\ ,\ \ 
H_{\alpha \beta \bar 0}(u) = h_{\alpha \beta}(u)\ .\tag5.1$$
\endproclaim
\medskip
Now
recall that for any $A \in \gl$ and any $u\in \LC$, 
the corresponding  fundamental vector field $A^*_u$ at $u$ is equal 
to the tangent vector at $t = 0$ of the curve $a(t) = u \circ e^{tA}$.
Therefore, by Lemma 3.6, 
the element $A$ belongs to the algebraic vertical subspace $\g_u$ 
if and only if 
$$\left.\frac{d}{dt}\left[\bh_{e^{tA}\cdot f_0}(e^{tA}\cdot e_\alpha, 
e^{tA} \cdot e_{\bar \beta})\right]\right|_{t=0} = 0\ .\tag5.2$$
Representing $A$ with the associated 
matrix  $A^\alpha_\beta$, condition (5.2) can be written as:
$$A^\gamma_\alpha \delta_{\gamma\beta} + 
\delta_{\alpha \gamma} A^{\bar \gamma}_{\bar \beta} + 
A^\gamma_0 H_{\alpha \bar \beta \gamma}(u) + 
A^{\bar \gamma}_{\bar 0}H_{\alpha \bar \beta \bar \gamma}(u) = 0\ .
\tag5.3$$
Using Lemma 5.1, we immediately obtain the 
following.\par
\proclaim{Proposition 5.2} For any $u\in \UF$, the algebraic
vertical subspace $\g_u \subset \goth{gl}_n(\C)$ is defined by 
the following  equations:
$$
\matrix
A^0_0 + A^{\bar 0}_{\bar 0} = 0\ , \qquad
A^0_\lambda + A^{\bar \lambda}_{\bar 0} + A^\nu_0 h_{\lambda \nu}(u) = 0\ ,\\
\phantom{AAA}\\
A^\mu_\lambda + 
 A^{\bar \lambda}_{\bar \mu} + 
A^\nu_0H_{\lambda \bar \mu \nu}(u) + 
A^{\bar \nu}_{\bar 0}H_{\lambda \bar \mu\bar\nu}(u) = 0\ .
\endmatrix\tag5.4$$
\endproclaim
Equations (5.4)  are called {\it the defining equations of the 
algebraic vertical subspace $\g_u$\/}.\par
\medskip
 For any $u \in \UF$, 
consider the following basis for $\g_u$ ($\lambda \geq \mu$):\par
$$t \= J_o \cdot E^0_0 \ ,\tag5.5$$
$$E^R_{\lambda, \mu} \= E^\lambda_\mu - E^\mu_\lambda 
\ ,
\quad E^I_{\lambda, \mu} \= J_o \cdot(E^\lambda_\mu + E^\mu_\lambda)\ ,\tag5.6$$
$$ E_{2\lambda}(u) =   E^0_\lambda - 
\left[\re\{h_{\lambda \mu}(u)\} + \delta_{\lambda \mu}\right] E^\mu_0 -
\im\{h_{\lambda \mu}(u)\} J_o \cdot E^\mu_0 -$$
$$ 
- \re\{H_{\bar \nu\mu \lambda }(u)\} E^\mu_\nu - 
\im\{H_{\bar\nu \mu \lambda }(u)\} J_o \cdot E^\mu_\nu \ ,\tag5.7$$
$$E_{2\lambda+1}(u) = J_o\cdot E^0_\lambda  - 
\left[\re\{h_{\lambda \mu}(u)\} - \delta_{\lambda \mu}\right]J_o\cdot E^\mu_0 +
\im\{h_{\lambda \mu}(u)\}  E^\mu_0 -$$
$$
- \re\{H_{\bar\nu \mu \lambda }(u)\} J_o\cdot E^\mu_\nu +
\im\{H_{\bar \nu \mu \lambda }(u)\} E^\mu_\nu \ .\tag 5.8$$
Consider also the  complex valued 
vector fields  
$\tilde e'_\lambda$, $\tilde e'_{\bar \lambda}$ defined as
$$
\tilde e'_\lambda|_u \=  \frac{1}{2}\left\{
E_{2\lambda}(u)^*|_u - \unim   E_{2\lambda+1}(u)^*|_u\right\} \ ,\qquad 
 \tilde e'_{\bar\lambda}|_u \= \overline{\tilde e'_\lambda|_u}\ .\tag5.9$$
If we consider the  vector fields $\tilde e'_\lambda$ and
the generalized fundamental vector fields 
$\tilde E^R_{\lambda, \mu}$,  
$\tilde E^I_{\lambda, \mu}$, $\tilde t$
as vector fields of $T^\C\LC$,
we may  write them as linear 
combinations of the vector fields  $\E^\alpha_\beta$ and
$\E^{\bar \alpha}_{\bar \beta}$. In this way, we obtain the following expressions:
$$
\tilde E^R_{\lambda, \mu} = (\E^\lambda_\mu - 
\E^\mu_\lambda) + (\E^{\bar \lambda}_{\bar \mu} - \E^{\bar \mu}_{\bar \lambda}) 
\ ,\qquad 
\tilde E^I_{\lambda, \mu} = \unim(\E^\lambda_\mu + 
\E^\mu_\lambda)   - 
\unim (\E^{\bar \lambda}_{\bar \mu} + \E^{\bar \mu}_{\bar \lambda}) 
\tag 5.10$$
$$
\tilde t = \unim \E^0_0 - \unim \E^{\bar 0}_{\bar 0}\ ,\qquad
\tilde e'_\lambda =
\E^0_\lambda|_u  - \E^{\bar \lambda}_{\bar 0}|_u -
 h_{\mu\lambda}(u) \E^\mu_0|_u   -
H_{\bar \nu \mu\lambda}(u) \E^\mu_\nu|_u\ .\tag 5.11$$
Notice that the vectors $E^R_{\lambda, \mu}$,
$E^I_{\lambda, \mu}$ constitute a basis for $\goth u_{n-1}$; 
hence
 any complex vector 
$X \in T^\C_u\Bbb V_x$, which is vertical
w.r.t. the 
projection $\hat \pi_{\pi(x)}:\Bbb V_{\pi(x)} \to S_{\pi(x)}
 = \Bbb V_{\pi(x)}/U_{n-1}$, is linear combination
of the vectors $\tilde E^R_{\lambda, \mu}(u)$,
$\tilde E^I_{\lambda, \mu}(u)$  and hence, by  (5.10), it is of the form
$$X = C^\mu_\nu \left(\E^\nu_\mu|_u - \E^{\bar \mu}_{\bar \nu}|_u\right)
\tag 5.12$$
for some uniquely determined $C^\mu_\nu \in \C$.\par
At the same time, for any $u\in \UF$,
the generalized fundamental vector fields 
$\tilde t$ and $\widetilde{E_{i}(u)}$, $ i = 2, \dots, 2n-1$, 
span a
subspace of $T_u \Bbb V_{\pi(u)}$, which is of dimension $2n-1$ and which is 
complementary  to the vertical distribution. More precisely we have the following:
\medskip
\proclaim{Lemma 5.3} For any $u\in \UF$
$$\tilde \epsilon_1|_u = \tilde t|_u\ \mod
 \text{span}_\R\{ A^*_u, \ A\in \goth u_{n-1}\}\ ,$$
$$\tilde \epsilon_{2\lambda} |_u = (E_{2 \lambda}(u))^*|_u  
\  \mod \text{span}_\R\{ A^*_u, \ A\in \goth u_{n-1}\}\ ,$$
$$\tilde \epsilon_{2\lambda+1} |_u = (E_{2 \lambda +1}(u))^*|_u  
\  \mod \text{span}_\R\{ A^*_u, \ A\in \goth u_{n-1}\}\ .
\tag5.13$$
Moreover,  
there exist some complex  valued functions 
$S^{\rho}_{\sigma}$ and 
$S^{\rho}_{\sigma \lambda}$ 
 such that 
$$\tilde \epsilon_1|_u =  \tilde t|_u -  S^{\rho}_{\sigma}
(u) (\E^\sigma_\rho|_u  - \E^{\bar\rho}_{\bar \sigma}|_u)
\ ,\tag 5.14$$
$$\tilde e_{\lambda}|_u =  \E^0_\lambda|_u - \E^{\bar \lambda}_{\bar 0}|_u 
-  h_{\sigma \lambda}(u)\E^\sigma_0|_u  -  H_{\bar \rho \sigma \lambda}(u)
\E^\sigma_\rho|_u -
S^{\rho}_{\sigma \lambda}(u)\left(\E^\sigma_\rho|_u -
\E^{\bar\rho}_{\bar \sigma}|_u\right)
\ .\tag5.15$$
\endproclaim
\demo{Proof} In order to prove  the claim, observe that, 
for any fiber $\Bbb V_x$ and 
 any  frame $u = \{f_0, \dots, f_{2n-1}\}\in \Bbb V_x$, 
the element $y = [u] = \hat \pi_x(u) \in  \Bbb V_x/U_{n-1} = S^{2n-1}$
may be identified with the first element  $f_0 = u(\epsilon_0)$ of the 
frame. Therefore 
the vector
$t^*|_u$
is projected by $\hat \pi_{x}$ onto the  vector
of
$T_{[u] = f_0} S_x$ given by 
$$\hat \pi_x{}_*(t^*|_u) = \left.\frac{d}{ds} (\hat \pi_x(u \circ 
e^{sJ_o E^0_0}))\right|_{s=0} =\left.\frac{d}{ds} (u \circ 
e^{sJ_o E^0_0}(\epsilon_0))\right|_{s=0} = $$
$$ = u \circ J_o E^0_0(\epsilon_0) = u \circ J_o E^0_0(\varepsilon_0 + 
\varepsilon_{\bar 0}) = u(i \varepsilon_0 - i \varepsilon_{\bar 0}) = 
f_1\ .$$
By a similar argument one can check that 
$(E_{2 \lambda}(u))^*|_u$ and $(E_{2\lambda+1}(u))^*_u$
are mapped by $\hat \pi_{x}$ onto the  vectors
$f_{2\lambda}$, $f_{2\lambda + 1} = J f_{2\lambda}$, respectively. This
proves (5.13).\par
From (5.13) and (5.12), the formulae (5.14) and (5.15) 
follow. \qed
\enddemo
In the following, in order to have more symmetry in some formula, we will
often write $H_{\bar 0\sigma \lambda}(u)$ in place
$h_{\sigma\lambda}(u)$, since they coincide by Lemma 5.1.
\bigskip
\subsubhead 5.3 The structure functions of the absolute parallelism of $\UF$
\endsubsubhead
\medskip
Cartan-Sternberg theorem (see [St]; we recall the complete
statement - which is indeed quite long - in the Appendix) 
implies  that a complete set
of invariant functions for an  absolute parallelism $\{X_1, \dots, X_m\}$
is given by the structure functions  $c^i_{jk}$,   defined by 
$[X_j, X_k] = c^i_{jk} X_i$, and by their derivatives  $X_{i_1}(
\dots X_{i_p}(c^i_{jk})\dots)$, with  $p$ less
or equal to some finite order $r$. 
The order $r$ in general depends on the absolute parallelism, but, 
in case of real analytic data, there
exists an upper bound for $r$ which depends only on the dimension of
the manifold (this is a consequence of Cartan-K\"ahler theorem; see [BCG]).
\par
From this remarks and Proposition 4.6 (1),  we conclude 
that the structure functions 
$c^i_{jk}$ 
 and the derivatives $X_{i_1}(
\dots X_{i_p}(c^i_{jk})\dots)$ of the absolute parallelism (4.3) 
are 
a complete system of invariant functions for
the complex Finsler manifold $(M,J, F)$. \par
\medskip
In this section we want to describe these structure functions.\par
\medskip
The structure functions $c^i_{jk}$ corresponding to
 Lie brackets
of  two generalized fundamental vector fields 
$\tilde X, \tilde Y$, with $X$ and $Y$ in 
$\{ t, 
{E^R_{\lambda, \mu}}, E^I_{\lambda, \mu}\}$, 
are  computed by  the Lie
brackets in $\u$. In fact, $\tilde X$ and $\tilde Y$ are 
the fundamental vector fields in the usual sense
and hence 
$[\tilde X, \tilde Y] = \widetilde{[X,Y]} $. In particular,
{\it for those Lie brackets, the corresponding
 structure functions $c^i_{jk}$
are   the structure constants
of the Lie algebra $\u$\/}. \par
\smallskip
The structure functions $c^i_{jk}$ corresponding to
Lie brackets  between a
 fundamental vector fields 
$\tilde X$, with $X\in\{
{E^R_{\lambda, \mu}}, E^I_{\lambda, \mu}\}$ 
and a vector field $Y$ in the set
$\{\hat \epsilon_i, \tilde \epsilon_j\}$ can be evaluated  recalling 
 that $\tilde X$ is  a fundamental
vector field  associated to an element in 
$\goth u_{n-1}$ and that $Y$ belongs either to $\H$ or 
to $\W$, which are both invariant under the action of $U_{n-1}$. This implies
that the action of $\tilde X$ on the set  
$\{\hat \epsilon_i, \tilde \epsilon_j\}$
is  equal to
the standard action of $X \in \goth u_{n-1}$ on the basis
$\{\epsilon_i, \epsilon_j\}$ of 
$V \oplus W = \C^n \oplus \C^{n-1}$ (see e.g. Prop. 2.3 in [KN], vol. I).
 Therefore, {\it also for these 
Lie brackets,  the corresponding
structure functions  are constant\/}. \par
\smallskip
The structure functions  corresponding to the Lie brackets
of  two vector fields in the set $\{\hat \epsilon_i, \tilde \epsilon_j, 
\tilde t\}$ are given by the real and imaginary 
parts of the Lie brackets described in the following Propositions 5.4,  
5.5 and 5.6. \medskip
\proclaim{Proposition 5.4} 
A Lie bracket of a pair of vector fields  in the set
$\{ \hat e_\beta, \hat e_{\bar \gamma}, \tilde t\}$ or of a pair in the set 
$\{
\tilde e_\lambda, \tilde e_{\bar \lambda}, 
\tilde t\}$ has one of the following forms:
$$[\hat e_\beta, \hat e_\gamma] =
-  T^\alpha_{\beta \gamma}\hat e_\alpha\ ,\tag5.16$$
$$[\hat e_\beta, \hat e_{\bar \gamma}] =
 -  R^\lambda_{0 \beta \bar \gamma}  \tilde e_{\lambda} +
R^0_{\lambda \beta \bar \gamma} \tilde e_{\bar \lambda} +
\unim   R^0_{0 \beta \bar \gamma} \tilde t
- R^\rho_{\sigma \beta\bar \gamma}(
\E^\sigma_\rho - 
\E^{\bar\sigma}_{\bar\rho})\ ,
\tag5.17
$$
$$
[\tilde t, \hat e_\beta] =  \unim   \delta_{\beta 0}\hat e_0\ ,\ 
\ [\tilde t, \hat e_{\bar\beta}] =  - \unim  \delta_{\beta 0}\hat e_{\bar 0}\ ,
\tag5.18$$
$$
[\tilde e_\mu, \tilde e_{\nu}] = 0 
\ ,\qquad
[\tilde e_\mu, \tilde e_{\bar\nu}] = - 
 \unim  \delta_{\mu \nu} \tilde t - (\delta^\rho_\mu \delta_{\sigma \nu} -
Q^{\rho}_{\sigma \mu \bar \nu})(\E^\sigma_\rho -
\E^{\bar\rho}_{\bar\sigma})
\ ,
\tag5.19
$$
$$
[\tilde t, \tilde e_\nu]_u =  - \unim   \tilde e_\nu - \unim  
Q^\rho_{\sigma 0 \nu} (\E^{\sigma}_\rho - \E^{\bar\rho}_{\bar\sigma}) \ ,\tag5.20$$
$$
[\tilde t, \tilde e_{\bar\nu}] = \unim  \tilde e_{\bar\nu}
-\unim   \overline{ Q^{\sigma}_{\rho 0 \nu}} (\E^{\sigma}_\rho 
-\E^{\bar\rho}_{\bar\sigma})\ ,
\tag5.21
$$
where  
 $T^\alpha_{\beta \gamma}$, $R^\alpha_{\beta \gamma \bar \delta}$,
$Q^{\rho}_{\sigma \mu \bar \nu}$, $Q^{\rho}_{\sigma 0 \nu}$
and  $Q^{\bar \rho}_{\bar \sigma 0 \nu}$
are some uniquely determined  $\C$-valued functions.
\endproclaim 
\demo{Proof}
Recall that the vector
fields $\hat e_\alpha$ are holomorphic vector fields
in  $\H^{10} \subset \H^\C$. Since the CR structure
$(\H, \hat J)$ is integrable,  
$[\hat e_\beta, \hat e_\gamma]$ takes
values in $\H^{10}$. From this
$(5.16)$ follows. \par
To prove (5.17) notice that 
 since $\pi: \LC \to M$ is 
holomorphic, $\pi_*(\hat e_\beta)$ and that
$\pi_*(\hat e_{\bar\gamma})$ are holomorphic and antiholomorphic in $(M,J)$, 
respectively.
Hence  $\pi_*([\hat e_\beta,\hat e_{\bar\gamma}]) = 
[\pi_*(\hat e_\beta),
\pi_*(\hat e_{\bar\gamma})] = 0$ and this proves that for any $u$, the vector
$[\hat e_\beta,\hat e_{\bar \gamma}]_u$ is a complex vertical vector. Since
at any $u$, the vectors 
$\tilde e_\lambda|_u$, $\tilde e_{\bar\lambda}|_u$, 
$\unim  \tilde t_u$ and $\E^\sigma_\rho|_u - \E^{\bar\rho}_{\bar\sigma}|_u$
are linearly independent over $\C$ and the complexified 
vertical subspace $\Cal V^\C_u \subset T^\C_u\UF$ is equal to their span,
(5.17) follows.
\par
To check  (5.18),  
recall that $\tilde t$
is the fundamental vector field in  $\LC$ associated to 
$J_o\cdot E^0_0$;  then the formula   follows  from  definitions
and the fact that $\H$ is $U_{n-1} \times T^1$-invariant. \par
(5.19) is a consequence of the properties of the distribution 
$\W$ (see Definition 4.3), of (5.14) and (5.15).\par
The proof of (5.20) is the following. Pick a frame $u_o \in \UF$ 
and let $X_\nu$
the complex vector field  in $T^\C L^\C(M)$ defined by 
$$X_\nu = \E^0_\nu - \E^{\bar \nu}_{\bar 0} 
- H_{\bar \alpha \sigma \nu}(u_o)\E^\sigma_\alpha  - 
S^{\rho}_{\sigma \lambda}(u_o)(\E^\sigma_\rho -
\E^{\bar\rho}_{\bar \sigma})\ .$$
From definitions and Lemma 5.3, one can check that
$$[\tilde t, \tilde e_\nu]_{u_o} = 
[\tilde t, X_\nu]_{u_o}  \ \mod span_\C \{\E^\sigma_0|_{u_o},
\E^\sigma_\rho|_{u_o},
\E^{\bar\sigma}_{\bar\rho}|_{u_o}\} \ .$$
By  the properties 
of the Lie brackets between fundamental vector fields in $\LC$, it follows that
there exists some complex functions $Q^A_{B C D}$ such that 
$$[\tilde t, \tilde e_\nu] = - \unim \tilde e_\nu
- \unim Q^0_{\sigma 0 \nu} \E^\sigma_0- \unim  
Q^\rho_{\sigma 0 \nu}  \E^{\sigma}_\rho -\unim   
Q^{\bar\rho}_{\bar \sigma 0 \nu} \E^{\bar\sigma}_{\bar\rho}\ .$$
Notice that by (5.10) and (5.11) the vector
$$[\tilde t, \tilde e_\nu]_{u} + \unim \tilde e_\nu|_{u} = 
- \unim Q^0_{\sigma 0 \nu}(u) \E^\sigma_0|_{u} - \unim  
Q^\rho_{\sigma 0 \nu}(u)\E^{\sigma}_\rho -\unim|_{u}  
Q^{\bar\rho}_{\bar \sigma 0 \nu}(u) \E^{\bar\sigma}_{\bar\rho}|_{u}$$
belongs to $span_\C\{\tilde e_\lambda|_{u}, \tilde e_{\bar \lambda}|_{u}, 
\tilde t|_u, \E^{\sigma}_\rho|_{u} - \E^{\bar\rho}_{\bar\sigma}|_{u}\}$ only if 
$ Q^0_{\sigma 0 \nu}(u) = 0$ and $Q^{\bar\rho}_{\bar \sigma 0 \nu}(u)
= Q^\sigma_{\rho 0 \nu}(u)$. This concludes the proof of
(5.20). 
(5.21) is obtained from (5.20) by conjugation and taking into account
 the fact that
$\overline{\tilde t} = \tilde t$.\qed
\enddemo
\medskip
\proclaim{Proposition 5.5} At any frame $u \in \UF$
and for any 
$\rho$, $\sigma$, $\lambda$, $\mu = 1, \dots ,n-1$, 
\roster
\item $S^{\rho}_{\sigma \lambda}(u) = 0\ ;$
in particular $\tilde e_\lambda' \equiv \tilde e_\lambda\ ;$
\item $Q^\rho_{\sigma 0 \mu}(u) =  0$ and 
$$Q^\rho_{\sigma \lambda \bar \mu}(u) = 
\tilde e_{\bar \mu}( H_{\sigma \lambda \bar \rho})|_u = 
H_{\bar \mu \bar \rho \sigma \lambda}(u)
- h_{\bar \mu \bar \rho}(u) h_{\sigma \lambda}(u)  - H_{\nu \bar \mu \bar \rho }(u)
H_{\sigma \lambda \bar \nu}(u) \ ;$$
\item the Lie brackets
  $[\tilde e_\lambda, \hat e_0]$,
$[\tilde e_{\bar \lambda}, \hat e_{0}]$,
$[\tilde e_\lambda, \hat e_\mu]$,
$[\tilde e_{\bar \lambda}, \hat e_{\mu}]$ have the following forms:
$$
[\tilde e_\lambda, \hat e_0] =  
\hat e_\lambda
 +\unim   P^{0}_{0 \lambda 0}\tilde t -
P^{\mu}_{0 \lambda 0} \tilde e_\mu +
P^{0}_{\mu \lambda 0}\tilde e_{\bar\mu}
-  P^{\sigma}_{ \rho \lambda 0} (\E^\rho_\sigma -
\E^{\bar\sigma}_{\bar\rho})\ ,
\tag5.22$$
$$
[\tilde e_\lambda, \hat e_\mu] = - H_{\bar \alpha \mu \lambda} \hat e_\alpha
+ \unim   P^{0}_{0 \lambda \mu} \tilde t -
P^{\nu}_{0 \lambda \mu} \tilde e_\nu +
P^{0}_{\nu \lambda \mu} \tilde e_{\bar\nu} 
-  P^{\sigma}_{ \rho \lambda \mu} (\E^\rho_\sigma -  
\E^{\bar\sigma}_{\bar\rho})
\ ,\tag5.23$$
$$[\tilde e_{\bar\lambda}, \hat e_0] =  
- \unim   P^{0}_{0 \bar\lambda 0} \tilde t +
P^{\mu}_{0 \bar\lambda 0} \tilde e_\mu -
P^{0}_{\mu \bar\lambda 0} \tilde e_{\bar\mu} 
+  P^{\sigma}_{ \rho \bar\lambda 0} (\E^\rho_\sigma - 
 \E^{\bar\sigma}_{\bar\rho})
\ ,
\tag5.24$$
$$
[\tilde e_{\bar \lambda}, \hat e_{\mu}] =  
- \delta_{\lambda \mu} 
\hat e_{0} 
 -\unim   P^{0}_{0 \bar\lambda \mu} \tilde t +
P^{\nu}_{0 \bar \lambda \mu} \tilde e_\nu -
P^{0}_{\nu \bar \lambda \mu} \tilde e_{\bar\nu} 
+ P^{\sigma}_{ \rho \bar\lambda  \mu} (\E^\rho_\sigma -  
 \E^{\bar\sigma}_{\bar\rho})
\ ,\tag5.25$$
where   
$P^{A}_{BCD}$ 
are some uniquely determined 
complex valued functions. 
\endroster
\endproclaim
\demo{Proof} (1) For any frame $u$, let us denote the vertical
subspace of $\UF$  by  $\Cal V_u$.  We claim that  
$$
[\tilde e_\lambda, \hat e_0]_u =  
\hat e_\lambda \mod \Cal V^\C_u\ ,
\ \ 
[\tilde e_\lambda, \hat e_\mu]_u = 
- (H_{\bar \alpha \mu \lambda}
+ S^{\alpha}_{\mu \lambda}) \hat e_\alpha   \mod  \Cal V^\C_u\ ,\tag5.26$$
$$[\tilde e_{\bar\lambda}, \hat e_0]_u = 0 \mod  \Cal V^\C_u\ ,\ \ 
[\tilde e_{\bar \lambda}, \hat e_{\mu}] =  
- \delta_{\lambda \mu} 
\hat e_{0} 
 +\overline{S^{\mu}_{\nu \lambda}}\hat e_{\nu}  
\mod  \Cal V^\C_u\ , \tag5.27$$
where $S^\sigma_{\rho \lambda}$ and $S^{\bar \sigma}_{\bar \rho \lambda}$
are the complex functions defined in Lemma 5.3 and we let
$S^0_{\rho \lambda} = 0$.\par
To prove (5.26), let us fix a frame $u_o$
and consider the complex vector 
field $X_\lambda$  in $T^\C L^\C(M)$ defined by 
$$X_\lambda = \E^0_\lambda - \E^{\bar \lambda}_{\bar 0} - 
(H_{\bar \alpha \sigma \lambda}(u_o) + S^{\alpha}_{\sigma \lambda}(u_o))
\E^\sigma_\alpha + 
S^{\sigma}_{\rho \lambda}(u_o)   
\E^{\bar\sigma}_{\bar \rho}\ .$$
Let us also extend  $\hat e_0$ to a  vector field  
on a neighborhood $\Cal U\subset L^\C(M)$ of $u_o$.\par
From (5.15) and from definitions we get
$$[\tilde e_\lambda, \hat e_0]_{u_o} = 
[X_\lambda, \hat e_0]_{u_o} \mod \text{span}_\C \{\E^\sigma_\alpha|_{u_o}, 
\E^{\bar\sigma}_{\bar\rho}|_{u_o}\} \ .$$
In particular, 
$\theta_{u_o}([\tilde e_\lambda, \hat e_0]) = 
\theta_{u_o}([X_\lambda, \hat e_0])$. Moreover,
 $X_\lambda(\theta(\hat e_0))|_{u_o} = 0$: in fact  
$X_\lambda|_{u_o}= \tilde e_\lambda|_{u_o}$ and hence it is tangent to $\UF$;
on the other hand $\theta(\hat e_0) = e_0$ at all points of $\UF$. \par
Therefore
$$\theta_{u_o}([X_\lambda, \hat e_0]) = 
-(\Cal L_{X_\lambda}\theta)_{u_o}(\hat e_0)
 = -(\Cal L_{\E^0_\lambda - \E^{\bar \lambda}_{\bar 0}}\theta)_{u_o}(\hat e_0) + $$
$$ +
(H_{\bar \alpha \rho \lambda}(u_o) + S^{\alpha}_{\rho \lambda}(u_o))
 \left(\Cal L_{\E^\rho_{\alpha}}\theta\right)_{u_o}(\hat e_0) 
-  S^{\rho}_{\sigma \lambda}(u_o) (\Cal L_{\E^{\bar\rho}_{\bar\sigma}}
\theta)_{u_o}(\hat e_0) =  e_\lambda \ .\tag5.28$$
This implies that  at all points the vectors $[\tilde e_\lambda, \hat e_0]$
and $\hat e_\lambda$ differ by a complex vertical vector. This implies  
the first identity in (5.26). The second identity in (5.26) and 
the two identities of (5.27) are proved with the same  arguments.\par
Now, we can prove 
that $S^{\lambda}_{\mu \nu} = 0$.  Indeed
from the Jacobi identities, Proposition 5.4, (5.26) and (5.27)
$$0 =\theta^\mu([\tilde e_\lambda,[\tilde e_{\bar \nu}, \hat e_{0}]]) + 
\theta^\mu([\tilde e_{\bar \nu},[\hat e_{0}, \tilde e_\lambda]]) + 
\theta^\mu([\hat e_{0},[\tilde e_\lambda, \tilde e_{\bar \nu}]])  = $$
$$ = \theta^\mu([\tilde e_{\bar \nu},[\hat e_{0}, \tilde e_\lambda]]) = 
 - \overline{S^{\lambda}_{\mu \nu}}\ .$$
\medskip
(2) is a consequence of (1), (5.26), (5.27), the Jacobi identities and  Lemma 2.4.
In fact,
$$0 = \theta^\rho([\tilde e_\lambda,[\tilde e_{\bar \mu}, \hat e_{\sigma}]]) + 
\theta^\rho([\tilde e_{\bar \mu},[\hat e_{\sigma}, \tilde e_\lambda]]) + 
\theta^\rho([\hat e_{\sigma},[\tilde e_\lambda, \tilde e_{\bar \mu}]])  = $$ 
$$ = \theta^\rho([\tilde e_\lambda,[\tilde e_{\bar \mu}, \hat e_{\sigma}]]) +
\theta^\rho([\tilde e_{\bar \mu},[\hat e_{\sigma}, \tilde e_\lambda]]) +  
\delta^\rho_\lambda \delta_{\sigma \mu} - Q^{\rho}_{\sigma \lambda \bar \mu}= $$
$$ =
- \delta_{\mu \sigma} \delta^\rho_\lambda + \tilde e_{\bar \mu}(
H_{\bar \rho \sigma\lambda}) 
+  
\delta^\rho_\lambda \delta_{\sigma \mu}  -  Q^{\rho}_{\sigma \lambda \bar \mu}$$
and 
$$0 = \theta^\rho([\tilde t,[\tilde e_{\bar \mu}, \hat e_{\sigma}]]) + 
\theta^\rho([\tilde e_{\bar \mu},[\hat e_{\sigma}, \tilde t]]) + 
\theta^\rho([\hat e_{\sigma},[\tilde t, \tilde e_{\bar \mu}]])  =  
\unim \ \overline{Q^{\sigma}_{\rho 0  \mu}}\ .$$
Now, using Lemma 2.4, (5.15) and 
the fact that $S^\rho_{\sigma \lambda} = 0$,
a straightforward computation 
shows  that
$$\tilde e_{\bar \mu}(
H_{\bar \rho \sigma\lambda}) = -
h_{\bar \mu \bar \rho} h_{\sigma \lambda}  - H_{\nu \bar \mu \bar \rho }
H_{\sigma \lambda \bar \nu} 
+ H_{ \bar \mu\bar \rho \sigma \lambda } $$
and this concludes the proof.\par
\medskip
(3) is an immediate consequence of (5.26), (5.27) and of claim (1).\qed
\enddemo
\medskip
\proclaim{Proposition 5.6} The structure functions 
$P^A_{BCD}$ 
defined in Proposition  5.5 are the following:
\roster
\item $P^{0}_{0 \rho \gamma} = 
P^{0}_{0 \bar \rho \gamma} = 
P^{0}_{\lambda \rho \gamma} = 
P^{0}_{\lambda \bar \rho \gamma} = 
P^\lambda_{0 \rho \gamma} = P^{\lambda}_{\mu \rho \gamma} = 0 $;
\item $P^\lambda_{0 \bar \rho \gamma} = - \hat e_\gamma (h_{\bar \lambda 
\bar \rho})$, $P^{\lambda}_{\mu \bar \rho \gamma} = 
- \hat e_{\gamma}(H_{\bar \lambda \mu \bar \rho})$.
\endroster
\endproclaim
\demo{Proof} It suffices to use the Jacobi 
identities, Proposition 5.4 and 
Proposition 5.5. In fact
$$0 = \theta^{\bar 0}([\tilde e_{\rho},[\hat e_{\bar 0}, \hat e_{\gamma}]]) + 
\theta^{\bar 0}([\hat e_{\gamma},[\tilde e_{\rho}, \hat e_{\bar 0}]]) + 
\theta^{\bar 0}([\hat e_{\bar 0},[\hat e_{\gamma}, \tilde e_{\rho}]])  = 
P^{0}_{0 \rho \gamma}
\ , $$
$$0 = \theta^{\bar 0}([\tilde e_{\bar \rho},[\hat e_{\bar 0}, \hat e_{\gamma}]]) + 
\theta^{\bar 0}([\hat e_{\gamma},[\tilde e_{\bar \rho}, \hat e_{\bar 0}]]) + 
\theta^{\bar 0}([\hat e_{\bar 0},[\hat e_{\gamma}, \tilde e_{\bar \rho}]])  = 
- P^{0}_{0 \bar \rho \gamma}
\ ,$$
$$0 = \theta^{\bar \lambda}([\tilde e_{\rho},[\hat e_{\bar 0}, \hat e_{\gamma}]]) + 
\theta^{\bar \lambda}([\hat e_{\gamma},[\tilde e_{\rho}, \hat e_{\bar 0}]]) + 
\theta^{\bar \lambda}([\hat e_{\bar 0},[\hat e_{\gamma}, \tilde e_{\rho}]])  = 
   P^{0}_{\lambda \rho \gamma} \ , $$
$$0 = \theta^{\bar \lambda}([\tilde e_{\bar \rho},[\hat e_{\bar 0}, \hat e_{\gamma}]]) + 
\theta^{\bar \lambda}([\hat e_{\gamma},[\tilde e_{\bar \rho}, \hat e_{\bar 0}]]) + 
\theta^{\bar \lambda}([\hat e_{\bar 0},[\hat e_{\gamma}, \tilde e_{\bar \rho}]])  = 
- P^{0}_{\lambda \bar \rho \gamma} \  ,$$
$$ 0 = \theta^{\bar 0}([\tilde e_{\rho},[\hat e_{\bar \lambda}, \hat e_{\gamma}]]) + 
 \theta^{\bar 0}([\hat e_{\gamma},[\tilde e_{\rho}, \hat e_{\bar \lambda}]]) +
 \theta^{\bar 0}([\hat e_{\bar \lambda},[\hat e_{\gamma}, \tilde e_{\rho}]]) =$$
$$ =
  P^{\lambda}_{0 \rho \gamma} - \overline{H_{\bar 0\lambda \sigma}} 
P^0_{\sigma \rho \gamma} =   P^{\lambda}_{0 \rho \gamma}\ ,$$
$$ 0 = \theta^{\bar \mu}([\tilde e_{\rho},[\hat e_{\bar \lambda}, \hat e_{\gamma}]]) + 
 \theta^{\bar \mu}([\hat e_{\gamma},[\tilde e_{\rho}, \hat e_{\bar \lambda}]]) +
 \theta^{\bar \mu}([\hat e_{\bar \lambda},[\hat e_{\gamma}, \tilde e_{\rho}]]) =
- P^{\lambda}_{\mu \rho \gamma}\ ,$$
$$ 0 = 
\theta^{\bar 0}([\tilde e_{\bar \rho},[\hat e_{\bar \lambda}, \hat e_{\gamma}]]) + 
 \theta^{\bar 0}([\hat e_{\gamma},[\tilde e_{\bar \rho}, \hat e_{\bar \lambda}]]) +
 \theta^{\bar 0}([\hat e_{\bar \lambda},[\hat e_{\gamma}, \tilde e_{\bar \rho}]]) =$$
$$ = - \hat e_{\gamma}
\left(\overline{H_{\bar 0 \lambda \rho}}\right)
- P^{\lambda}_{0 \bar \rho \gamma} 
 =  - \hat e_{\gamma}
(h_{\bar \lambda \bar \rho})
- P^{\lambda}_{0 \bar \rho \gamma}\ ,$$
$$ 0 = \theta^{\bar \mu}([\tilde e_{\bar\rho},[\hat e_{\bar \lambda}, \hat e_{\gamma}]]) + 
 \theta^{\bar \mu}([\hat e_{\gamma},[\tilde e_{\bar\rho}, \hat e_{\bar \lambda}]]) +
 \theta^{\bar \mu}([\hat e_{\bar \lambda},[\hat e_{\gamma}, \tilde e_{\bar\rho}]]) =$$
$$
= - \hat e_{\gamma}\left(
\overline{H_{\bar \mu \lambda\rho}}\right) -
P^{\lambda}_{\mu \bar \rho \gamma}   = 
- \hat e_{\gamma}(H_{\mu
\bar \lambda \bar \rho})  - P^{\lambda}_{\mu \bar \rho \gamma}  \ .\qquad
\qed$$
\enddemo
\bigskip
By the previous remarks and 
Propositions 5.4, 5.5 and 5.6, we now have the 
complete list for  the structure functions
of $\sigma^\H$ and  they  generate a complete 
system of invariant functions. We summarize the results
in the next corollary.
For notation and 
 indexing conventions, see \S 5.1.\par
\bigskip
\proclaim{Corollary 5.7} The structure functions $c^i_{jk}$ of the absolute
parallelism $\sigma^\H$, associated to the non-linear Hermitian connection
of a complex Finsler manifold $(M,J,F)$, are the following:
\roster 
\item"i)"  the structure constants of $(\goth u_{n-1} + \R) \ltimes \C^n$, where 
the action of  $(\goth u_{n-1} + \R)$  on $\C^n$ is the one induced by 
the standard action of $\goth u_n$ on $\C^n$;
\item"ii)" the real and the imaginary parts of the functions
$h_{\lambda \mu}$ and $H_{\lambda \mu \bar \nu}$;
\item"iii)" 
the real and the imaginary parts of  the functions 
$$P^\lambda_{0 \bar \rho \gamma} = - \hat e_\gamma (h_{\bar \lambda 
\bar \rho})\ ,\ \ 
P^\lambda_{\mu \bar \rho \gamma} = - \hat e_{\gamma}(H_{\bar \lambda \mu \bar \rho})\ ,
\ \ 
Q^\rho_{\sigma\lambda\bar \mu} =   H_{\bar \mu \bar \rho \sigma \lambda} 
- h_{\bar \mu \bar \rho} h_{\sigma \lambda} -
 H_{\nu \bar \mu \bar \rho }
H_{\sigma \lambda \bar \nu} 
\ ;$$
\item"iv)" the real and the imaginary parts
of the functions $T^\alpha_{\beta \gamma}$ and $R^\alpha_{\beta \gamma \bar \delta}$
defined by (5.12) and (5.13).
\endroster
\endproclaim
\bigskip
\bigskip
\subsubhead 5.4 The structure equations of a complex Finsler manifold
\endsubsubhead
\medskip
The structure equations of $\UF$ consist in the 
identities verified by the tautological 
1-form $\theta$, the connection form $\omega$ and  
the  differentials of their components.  They are direct consequences 
of 
the defining equations of the algebraic vertical subspaces
and  the structure functions of the absolute parallelism. 
The results are in 
the following
theorem.\par
For notation and 
 indexing conventions, see \S 5.1 and \S 5.2.\par
\medskip
\proclaim{Theorem 5.8}
Let $(M,J,F)$ be a complex Finsler space and let $\omega$  the connection
1-form on $\UF$ 
associated with  the non-linear Hermitian connection $\H$ of
$(M,J,F)$.
\roster
\item"i)" 
The holomorphic and anti-holomorphic components of  $\omega$ 
verify:
$$\omega^0_0 + \omega^{\bar 0}_{\bar 0} = 0\ ,\quad
\omega^0_\lambda +  \omega^{\bar \lambda}_{\bar 0} + h_{\lambda \nu} \omega^{\nu}_0
 = 0\ ,\quad
\omega^\lambda_\mu +
\omega^{\bar \mu}_{\bar \lambda} + 
H_{\bar \lambda\mu \nu} \omega^\nu_0 +
H_{\bar \lambda\mu\bar \nu} \omega^{\bar\nu}_{\bar 0} = 0\ .\tag5.29$$
\item"ii)" Let $\varpi^\alpha_\beta$ be 
the $\C$-valued 1-forms on $\UF$ given by
$$\varpi^0_0 = \omega^0_0\ ,\ \ 
\varpi^\lambda_0 = \omega^\lambda_0\ ,\ \ 
\varpi^0_\lambda = - \omega^{\bar \lambda}_{\bar 0}\ , 
\ \ 
\varpi^\mu_\nu = \omega^\mu_\nu + 
H_{\bar \mu \nu \lambda} \omega^{\lambda}_0\ ,\ \ 
\varpi^{\bar \alpha}_{\bar \beta} = \overline{\varpi^\alpha_\beta}\ .$$
Then $\varpi^\alpha_\beta$ verify:
$$\varpi^{\bar \alpha}_{\bar \beta} = - \varpi^\beta_\alpha\ ,\tag 5.30$$
$$\varpi^\lambda_0(\tilde e_\mu) = \delta_{\lambda \mu}\ ,
\ \ \varpi^\lambda_0(\tilde e_{\bar \mu}) = 0\ ,\quad
\varpi^0_\lambda(\tilde e_\mu) = 0\ ,\ 
\ \varpi^0_{\lambda}(\tilde e_{\bar \mu}) = - \delta_{\lambda \mu}\ ,\tag 5.31$$
$$\varpi^0_\lambda(\tilde t) = \varpi^0_\lambda(\tilde E^R_{\mu, \nu}) = 
\varpi^0_\lambda(\tilde E^I_{\mu, \nu}) = 0\ , 
\ \ \varpi^\lambda_0(\tilde t) = \varpi^\lambda_0(\tilde E^R_{\mu, \nu}) = 
\varpi^\lambda_0(\tilde E^I_{\mu, \nu}) = 0\ ;\tag 5.32$$
\item"iii)" The differentials of 
the tautological 1-form
$\theta$ and of the $\C$-valued 1-forms  
$\varpi^\beta_\gamma$ 
 are given by the following identities:
$$d\theta^\alpha + \varpi^\alpha_\beta \wedge \theta^\beta =
\Theta^\alpha + \Sigma^\alpha\ ;\tag5.33$$
$$d\varpi^0_0 + \varpi^0_\beta\wedge\varpi^\beta_0 = 
\Omega^0_0 \ ;\tag5.34$$
$$
d\varpi^\lambda_0 + 
\varpi^\lambda_\beta\wedge\varpi^\beta_0 =
\Omega^\lambda_0 + \Pi^\lambda_0\ 
,\quad
\ d\varpi^{0}_{\lambda} + 
\varpi^{0}_\beta\wedge\varpi^\beta_\lambda =
\Omega^0_\lambda + \Pi^0_\lambda\  ;
\tag5.35
$$
$$
d\varpi^\lambda_\mu + 
\varpi^\lambda_\beta\wedge\varpi^\beta_\mu = \Omega^\lambda_\mu
+ \Pi^\lambda_\mu + \Phi^\lambda_\mu\ ;
\tag5.36$$
where 
 $\Theta^\alpha$, $\Sigma^{\alpha}$, 
$\Omega^\alpha_\beta$, $\Pi^\lambda_0$,
$\Pi^0_\mu$, $\Pi^\lambda_\mu$ and 
$\Phi^\lambda_\mu$
are the the following $\C$-valued 2-forms:
$$\Theta^\alpha = 
\frac{1}{2}
T^\alpha_{\beta \gamma} \theta^\beta\wedge\theta^\gamma\ ,\quad
\Sigma^\alpha = 
H_{\bar \alpha \mu \lambda}\varpi^\lambda_0\wedge \theta^\mu \ ,
\quad
\Omega^\alpha_\beta = 
R^{\alpha}_{\beta \gamma \bar \delta}\theta^\gamma\wedge\theta^{\bar \delta}
\ ,\tag 5.37$$
$$ \Pi^0_\lambda =  -\hat e_{\bar \gamma}(h_{\lambda \rho})\varpi^\rho_0 \wedge 
\theta^{\bar \gamma} \ ,
\ \ \Pi^\lambda_0 =  -\hat e_\gamma(h_{\bar \lambda \bar \rho})
\varpi^0_\rho \wedge 
\theta^{\gamma} 
\ , \tag5.38$$
$$
\Pi^\lambda_\mu= 
- \hat e_\gamma(H_{\bar \lambda \mu \bar \rho}) \varpi^{0}_\rho \wedge
\theta^\gamma - 
\hat e_{\bar \gamma}(H_{\bar \lambda \mu \rho})
\varpi^{\rho}_{0} \wedge\theta^{\bar \gamma}\ ,\tag5.39$$
$$\Phi^\lambda_\mu = 
\left( H_{\bar \lambda \bar \sigma \mu \rho }
-  h_{\bar \lambda \bar \sigma} h_{\mu \rho}  - 
 H_{\nu \bar \lambda  \bar \sigma} H_{\bar \nu \mu \rho }
 \right)
\varpi^\rho_0 \wedge \varpi^0_\sigma
\ .\tag5.40$$
where $T^\alpha_{\beta \gamma}$ and $R^{\alpha}_{\beta \gamma \bar \delta}$
are the complex functions defined in (5.12) and (5.13).
\endroster
\endproclaim
We call the equations (5.29) and (5.33) - (5.36)  {\it the 
structure equations of the 
non-linear Hermitian connection of $(M,J,F)$\/}. 
\demo{Proof} 
 (5.29)  follows from the defining equations
of the algebraic vector subspaces. To check (5.30) - (5.32),  one  has only
to use the definitions, Lemma 5.3 and Proposition 5.5 (1). 
The structure equations (5.33) - (5.36) are proved by 
evaluating both sides  on all possible pairs 
of vector fields of the absolute parallelism (4.3) and 
checking that both sides give the same result. This 
can done directly by using Propositions 5.4, 5.5, 5.6
and formulae (5.31) and  (5.32). \qed
\enddemo
\bigskip
Consider  the following 2-forms on $\UF$ with values in $\gl \oplus \C^n$: 
$$\Theta = \sum_{\alpha = 0}^{n-1} e_\alpha\otimes \Theta^\alpha\ ,\qquad
\Sigma = \sum_{\alpha = 0}^{n-1} e_\alpha\otimes \Sigma^{\alpha}\ ,\qquad
\Omega = \sum_{\alpha, \beta = 0}^{n-1}E^\beta_\alpha\otimes \Omega^\alpha_\beta\ ,$$
$$
\Pi = \sum_{\lambda = 1}^{n-1}
E^0_\lambda\otimes \Pi^\lambda_0 +
\sum_{\mu = 1}^{n-1}
E^\mu_0\otimes \Pi^0_\mu + 
\sum_{\lambda,\mu = 1}^{n-1}
E^\mu_\lambda\otimes \Pi^\lambda_\mu\ ,\ \ 
\Phi = \sum_{\lambda, \mu = 1}^{n-1}
E^\mu_\lambda\otimes \Phi^\lambda_\mu\ $$
We  call $\Theta$  the
{\it (pure) torsion form\/}   and $\Sigma$ the {\it 
Finsler torsion form\/}. 
The 2-form $\Omega$ is called the {\it (pure) curvature form\/};
finally we name $\Pi$ and $\Phi$  
{\it oblique Finsler curvature\/} and  {\it vertical
Finsler curvature\/}, respectively.\par 
\medskip
The Finsler curvature and torsion forms
are non-zero only if the Finsler metric
is not associated with an Hermitian metric. The
following Proposition gives an important criterion to see
when this occurs.\par
\medskip
\proclaim{Proposition 5.9}  A complex Finsler metric $F$ 
is associated with an Hermitian metric $g$ if and only if the 
component $\Sigma^0$ of the
Finsler torsion form 
vanishes identically.\par
In this case,  $\Sigma = 0$, $\Pi = 0$ and $\Phi = 0$
 and $\Theta$ and $\Omega$
coincide with the torsion form and  the curvature forms of the linear 
Hermitian connection of $(M, J,g)$, respectively.
\endproclaim
\demo{Proof} From definitions, 
$\Sigma^0$ vanishes if and only if for any $x\in M$, any 
$0\neq v\in T_xM$ 
and any two trivially extended vector fields $X, Y \in T(T_xM)$,
$\bh(X^{10}, Y^{10}) = 0$. 
By Lemma 2.5, this occurs if and only if 
$F^2$ is associated
to an Hermitian metric $h$. The other part of the claim
follows immediately from the identity between $\UF$ and 
 the unitary frame bundle corresponding to the Hermitian metric associated
with $F$.\qed
\enddemo
Taking the exterior differential of both sides of the structure 
equations, one can obtain several identities that must be satisfied
by the structural functions and by  the torsion and curvature forms.
Some of them are given in the following Proposition.
When $F$ is associated with 
an Hermitian metric, they reduce to 
the usual symmetry identities and to the  Bianchi identities
of the torsion and curvature of a 
linear Hermitian connection.\par
\medskip
\proclaim{Proposition 5.10} Let $\Theta$ and $\Omega$ the pure 
curvature forms of a complex Finsler manifold $(M,J,F)$. 
Then the
components $\Theta^\alpha = T^\alpha_{\beta \gamma}
\theta^\beta\wedge\theta^\gamma$ and 
$\Omega^\alpha_\beta = R^\alpha_{\beta \gamma \bar \delta}
\theta^\gamma \wedge \theta^{\bar \delta}$ verify the following identities: \par
$$\Omega^\alpha_\beta = - \overline{\Omega^\beta_\alpha}\ ;\tag 5.40$$
\noindent{\rm (First Bianchi Identities)\/}
$$\hat e_\beta(T^\alpha_{\gamma \delta}) + 
\hat e_\gamma(T^\alpha_{\delta \beta}) +  
\hat e_\delta(T^\alpha_{\beta \gamma}) + 
 T^\alpha_{\varepsilon \beta} T^\varepsilon_{\gamma\delta} + 
T^\alpha_{\varepsilon \delta} T^\varepsilon_{\beta\gamma} +
 T^\alpha_{\varepsilon \gamma} T^\varepsilon_{\delta\beta} = 0\ ,\tag 5.41$$
$$R^{\alpha}_{\beta \gamma \bar \delta} - R^{\alpha}_{\gamma \beta \bar \delta}
- \hat e_{\bar \delta}(T^\alpha_{\beta \gamma}) -
H_{\bar \alpha \lambda \beta } R^\lambda_{0 \gamma \bar \delta} + 
H_{\bar \alpha \lambda \gamma } R^\lambda_{0 \beta \bar \delta} = 0\ ,\tag 5.42$$
\noindent {\rm (Second Bianchi Identities)\/}
$$\hat e_\gamma(R^\alpha_{\beta \delta \bar \varepsilon})
- \hat e_\delta(R^\alpha_{\beta \gamma \bar \varepsilon}) + 
R^\alpha_{\beta \zeta \bar \varepsilon} T^{\zeta}_{\gamma \delta} 
 + \hat e_{\gamma}(H_{\bar \alpha \beta \bar \lambda})
R^0_{\lambda \delta \bar \varepsilon} -
\hat e_{\delta}(H_{\bar \alpha \beta \bar\lambda})
R^0_{\lambda \gamma \bar \varepsilon}
= 0\ ,\tag 5.43$$
$$\hat e_{\bar \delta}(R^\alpha_{\beta \gamma \bar \varepsilon})
- \hat e_{\bar \varepsilon}(R^\alpha_{\beta \gamma \bar \delta}) -
R^\alpha_{\beta \gamma \bar \zeta} \overline{T^{\zeta}_{\delta
\varepsilon}}  - \hat e_{\bar \varepsilon}(H_{\bar \alpha \beta \lambda})
R^\lambda_{0 \gamma \bar \delta} + 
\hat e_{\bar \delta}(H_{\bar \alpha \beta \lambda})
R^\lambda_{0 \gamma \bar \varepsilon} = 0\ .\tag 5.44$$
\endproclaim
\demo{Proof} (1) The identity (5.30) implies that $d \varpi^\alpha_\beta = 
- \overline{d \varpi^\beta_\alpha}$. Then (5.40) follows directly
from this and 
the structure equations.\par
Now, consider the  exterior differential of both sides
of 
the structure equations (5.33) - (5.36). 
It can be easily checked that the 3-form
$$(\Omega^\alpha_\beta - H_{\bar \alpha \lambda \beta} \Omega^\lambda_0)\wedge
\theta^\beta - d\Theta^\alpha\ ,\quad
d\Omega^0_0\ ,\quad
d\Omega^\lambda_0 + d\Pi^\lambda_0\ ,\quad
d\Omega^0_\lambda + d\Pi^0_\lambda\ ,\quad
d\Omega^\lambda_\mu + d\Pi^\lambda_\mu$$
vanish identically on any three vectors of horizontal distribution $\H$. This 
implies that the components  of these 3-forms 
with respect to the forms
$$\theta^\beta\wedge \theta^\gamma\wedge \theta^\delta\ ,\  
\ \theta^\beta\wedge \theta^\gamma\wedge \theta^{\bar \delta}\ , 
\ \theta^\beta\wedge \theta^{\bar \gamma}\wedge \theta^{\bar \delta}\ ,
\ \theta^{\bar \beta}\wedge \theta^{\bar \gamma}\wedge \theta^{\bar \delta}$$
have to vanish identically. Such components are exactly the 
left hand sides of (5.41) - (5.44). \qed
\enddemo
\medskip
\remark{Remark 5.11} The exterior differentiation of 
the structure equations give 
several other identities for the torsion and the curvature, which are
not listed in Proposition 5.10.\par
When the Finsler metric is associated
with an Hermitian metric, those identities 
express only the property that
$T^\alpha_{\beta \gamma}(u)$ and 
$R^\alpha_{\beta \gamma \bar \delta}(u)$ depend on the frame $u = \{e_i\}$ 
as the components of some suitable tensor fields on $M$.
But when the Finsler metric $F$ is not associated with an 
Hermitian metric, these  identities  
give new (and somehow unexpected) relations between
$\Theta$, $\Omega$ and  
the  Finsler torsion and  curvatures $\Sigma$, $\Pi$ and $\Phi$.\par
One can  obtain a complete list of these identities (and 
avoid several tedious computations) using some symbolic
manipulation computer program.
\endremark
\bigskip
\bigskip
\subhead 6. The structure equations of  Lempert manifolds
\endsubhead
\bigskip
\subsubhead 6.1 Geodesics of a complex Finsler manifold
\endsubsubhead
\medskip
We continue to use all conventions given 
in \S 5.1, \S 5.2 and \S 5.4.\par
\smallskip
Let $\gamma: [a,b] \subset \R \to M$ be a smooth regular curve, 
that is a smooth curve 
such that the tangent vector $\dot \gamma_t$ is different
from $0$ for any $t\in [a,b]$. We call {\it lift of $\gamma$\/}
any smooth curve $\tilde \gamma: [a, b] \to \UF$ such that:
\roster 
\item"a)" $\pi\circ \tilde \gamma = \gamma$;
\item"b)" 
for any $t\in [a,b]$, the  
frame $\tilde \gamma_t = \{f_0(t), \dots, f_{2n-1}(t)\} \subset 
T^{10}_{\gamma_t}M$ is 
so that $f_0(t) \in \C^* \dot \gamma_t$. 
\endroster
Notice that 
a curve $\tilde \gamma: [a, b] \to \UF$
is a lift of $\gamma$ if and only if it projects onto $\gamma$
and there exists 
a smooth map $\varphi: [a,b] \to \R$ 
 such that for any frame $\tilde \gamma_t = \{f_i(t)\}$
$$\dot \gamma_t =  F(\dot \gamma_t) e^{\varphi_t J}\cdot \left(f_0(t) \right)= 
 F(\dot \gamma_t)\left(e^{i\varphi_t}e_0(t) + 
 e^{-i\varphi_t}e_{\bar 0}(t)\right)\ .\tag6.1$$
We call {\it  length of $\gamma$ \/} and {\it  energy of $\gamma$\/}
the integrals $L(\gamma)$ and
$E(\gamma)$ defined by
$$L(\gamma) = \int_a^b F(\dot \gamma_t) dt\ ,\quad
E(\gamma) = \int_a^b F^2(\dot \gamma_t) dt\ .\tag 6.2$$ 
Note that if $\tilde \gamma$ is any lift of $\gamma$ on $\UF$, then
$$L(\gamma) = \int_a^b \sqrt{\theta^0(\dot {\tilde \gamma}_t) 
\theta^{\bar 0}(\dot {\tilde \gamma}_t)}dt\ ,
\quad E(\gamma)  = \int_a^b \theta^0(\dot {\tilde \gamma}_t) 
\theta^{\bar 0}(\dot {\tilde \gamma}_t)dt\ .\tag 6.3$$
We recall that  a {\it variation with fixed endpoints\/}
of $\gamma$ is a smooth map $V: (-\delta, \delta) \times [a,b]
\to M$ such that 
\roster
\item $V(0,t) = \gamma_t$ for all $t\in [a,b]$;
\item for any $s\in (-\delta,\delta)$, the curve $\gamma^{(s)} = V(s, *)$
is a regular curve such that $\gamma^{(s)}_a = \gamma_a$ and 
$\gamma^{(s)}_b = \gamma_b$.
\endroster
\definition{Definition 6.1}
A regular curve $\gamma: [a, b] \to M$ is called {\it geodesic of the complex
Finsler manifold $(M,J,F)$\/}
if for any variation $V$ with fixed endpoints, the family 
of curves $\gamma^{(s)} = V(s, *)$
is so that 
$$\left.\frac{d E(\gamma^{(s)})}{dt}\right|_{s=0} = 0\ .\tag6.4$$
\enddefinition
The equations of Euler-Lagrange for a geodesic of a complex
Finsler manifold are given in the following Theorem.\par
\smallskip
\proclaim{Theorem 6.2} Let $\gamma: [a, b] \to M$ be a regular curve. 
\roster
\item For any   lift  $\tilde \gamma: [a, b] \to \UF$ of $\gamma$
consider the complex functions
$A(\tilde \gamma)_t$, $B(\tilde \gamma)_t$ and $C(\tilde \gamma)_t$
defined by 
$$A(\tilde \gamma)_t = \varpi^0_0(\dot{\tilde \gamma}_t)
\theta^{\bar 0}(\dot{\tilde \gamma}_t) - 
\left.\frac{d \theta^{\bar 0}
(\dot{\tilde \gamma}_s)}{ds}\right|_t\ ,$$
$$
B(\tilde \gamma)_t = 
\varpi^0_\lambda(\dot{\tilde \gamma}_t)
 + T^0_{\lambda 0}|_{\tilde \gamma_t} 
\theta^0
(\dot{\tilde \gamma}_t)\ ,\qquad
C(\tilde \gamma)_t = 
\varpi^{\bar 0}_{\bar \lambda}
(\dot{\tilde \gamma_t}) + \overline{T^{0}_{\lambda 
0}}|_{\tilde \gamma_t} 
\theta^{\bar 0}
(\dot{\tilde \gamma})\ .$$
They vanish identically for one lift of $\gamma$
if and only if they vanish for any lift
of $\gamma$.
\item $\gamma$ is a geodesic if and only if for any lift $\tilde \gamma$
and any  $t\in [a,b]$
$$\frac{d F(\dot \gamma)}{dt} = 0\ ,
\quad\varpi^0_0(\dot{\tilde \gamma_t})
\theta^{\bar 0}(\dot{\tilde \gamma}_t) - \left.\frac{d \theta^{\bar 0}
(\dot{\tilde \gamma}_s)}{ds}\right|_t  = 0\ ,\tag6.5$$
$$\varpi^0_\lambda(\dot{\tilde \gamma_t}) + T^0_{\lambda 0}|_{\tilde \gamma_t} 
\theta^0
(\dot{\tilde \gamma}_t) = 0\ ,\quad
\varpi^{\bar 0}_{\bar \lambda}
(\dot{\tilde \gamma_t}) + \overline{
T^{0}_{\lambda 0}}|_{\tilde \gamma_t} 
\theta^{\bar 0}
(\dot{\tilde \gamma}_t)= 0\ .\tag6.6$$
\endroster
\endproclaim
\demo{Proof} (1) Let $\tilde \gamma$ and $\tilde \gamma'$ be two lifts
of $\gamma$. Let also $\varphi_t$ and $\varphi'_t$  
two real functions so that (6.1) holds 
for $\tilde \gamma$ and $\tilde \gamma'$, respectively. 
This means that  $\tilde\gamma'_t=
 \tilde \gamma_t \circ \left(
e^{(\varphi_t - \varphi'_t)J} \right) $ and that 
$\dot{\tilde \gamma'}_t = \left(R_{e^{(\varphi_t - \varphi'_t)J}}\right)_*
\left(\dot{\tilde \gamma}_t\right)$ for any $t$. From the invariance
properties
of $\theta$, $\omega$ and of the torsion 2-forms $\Theta^\alpha$
 under the action 
of $T^1$, it follows immediately that 
$A(\tilde \gamma') = e^{\unim   (\varphi'_t - \varphi_t)} A(\tilde \gamma)$, 
$B(\tilde \gamma') = e^{\unim   (\varphi_t - \varphi'_t)}B(\tilde \gamma)$ and
$C(\tilde \gamma') = e^{\unim   (\varphi'_t - \varphi_t)}C(\tilde \gamma)$.\par
\smallskip
(2)  Consider a variation $V: (-\delta, \delta)\times[a,b] \to M$ 
of $\gamma$ with fixed endpoints and let 
 $\tilde V : (-\delta, \delta) \times [a,b]\to \UF$ be a smooth 
map  such that 
for any $s\in (-\delta, \delta)$, the curve $\tilde \gamma^{(s)} = 
\tilde V(s, *)$ is a lift of the curve $\gamma^{(s)} = V(s, *)$. \par
Let also $X$ and $Y$ the vector fields, which are tangent to $\tilde V
((-\delta, \delta) \times [a,b]) \subset \UF$, defined by 
$$X = \tilde V_*\left(\frac{\partial}{\partial t}\right)\ ,\quad
Y = \tilde V_*\left(\frac{\partial}{\partial t}\right)\ .$$
Note that $[X,Y] = \tilde V_*\left(\left[
\frac{\partial}{\partial t},\frac{\partial}{\partial s}\right]\right) = 0$.\par
From definitions and the commuting property of
$X$ and $Y$ we get that
$$\left.\frac{d E(\gamma^{(s)})}{ds}\right|_{s=0} = 
\int_a^b 
\left.Y\left(\theta^0(X) \theta^{\bar 0}(X)\right)\right|_{\tilde \gamma_t}dt = 
\int_a^b\left. \left
[d\theta^0(Y,X) \theta^{\bar 0}(X) + \right.\right.$$
$$ +
\theta^{0}(X) d\theta^{\bar 0}(Y,X)  
\left.\left.+X\left(\theta^0(Y)\right)
\theta^{\bar 0}(X) + \theta^{0}(X)X\left(\theta^{\bar 0}(Y)\right)\right]
\right|_{\tilde \gamma_t}dt = $$
$$ = \int_a^b\left. \left
[d\theta^0(Y,X) \theta^{\bar 0}(X) + 
\theta^{0}(X) d\theta^{\bar 0}(Y,X)  - 
\theta^{0}(Y)X\left(\theta^{\bar 0}(X)\right) -\right.\right.$$
$$\left.\left.
- \theta^{\bar 0}(Y)
X\left(\theta^0(X)\right)\right]\right|_{\tilde \gamma_t}dt + 
\int_a^b \left. X\left(\theta^{0}(Y)\theta^{\bar 0}(X) +
\theta^{0}(X)\theta^{\bar 0}(Y)\right)\right|_{\tilde \gamma_t}dt\ .\tag6.7$$
 Now,
$$\int_a^b \left. X\left(\theta^{0}(Y)\theta^{\bar 0}(X) +
\theta^{0}(X)\theta^{\bar 0}(Y)\right)\right|_{\tilde \gamma_t}dt = 
\left.\left(\theta^{0}(Y)\theta^{\bar 0}(X) +
\theta^{0}(X)\theta^{\bar 0}(Y)\right)\right|^{{\tilde \gamma}_b}
_{{\tilde \gamma}_a} = 0\tag6.8$$
because $\theta^0(Y)_{\tilde
\gamma_a} = \theta^0(Y)_{\tilde
\gamma_b} = 0$ since $V$ is a variation with fixed endpoints. \par
Using (6.8),  the fact that $\theta^\alpha(X) = 
\delta^\alpha_0 \theta^0(X)$ and the structure equations (5.29)
and (5.33),
we get
$$\left.\frac{d E(\gamma^{(s)})}{ds}\right|_{s=0} = 
\int_a^b\left
\{[- (\varpi^0_\alpha\wedge \theta^\alpha)(Y,X) + 
\Theta^0(Y,X) + \Sigma^0(Y,X)]\theta^{\bar 0}(X) + \right.$$
$$ + \theta^{ 0}(X)[- (\omega^{\bar 0}_{\bar \alpha}\wedge 
\theta^{\bar\alpha})(Y,X) + 
\Theta^{\bar 0}(Y,X) + \Sigma^{\bar 0}(Y,X)] - $$
$$ -\theta^{0}(Y)X (\theta^{\bar 0}(X))
\left.- \theta^{\bar 0}(Y)
X(\theta^0(X))\right\}_{\tilde \gamma_t}dt = $$
$$ = 
\int_a^b\left\{\theta^0(X) \theta^{\bar 0}(X)\right\}
_{\tilde \gamma_t}
\left\{\theta^0(Y) \left[\frac{\varpi^0_0(X)}{\theta^0(X)}
- \frac{X(\theta^{\bar 0}(X))}{
\theta^0(X) \theta^{\bar 0}(X)}\right] + \right.$$
$$
+ \theta^{\bar 0}(Y) \left[\frac{\varpi^{\bar 0}_{\bar 0}(X)}
{\theta^{\bar 0}(X)}
- \frac{X(\theta^{0}(X))}{\theta^0(X) \theta^{\bar 0}(X)}\right]+
\theta^\lambda(Y) \left[\frac{\varpi^0_\lambda(X)}{\theta^0(X)} 
+ T^0_{\lambda 0}\right] + $$
$$\left. +
\theta^{\bar \lambda}(Y) \left[\frac{\varpi^{\bar 0}_
{\bar \lambda}(X)}{\theta^{\bar 0}(X)}
+ T^{\bar 0}_{\bar \lambda\bar 0}\right]
\right\}_{\tilde \gamma_t}dt \tag6.9$$
Hence $\gamma$ is a geodesic if and only if   equation
(6.6) and the following equations (6.10) are identically satisfied:
$$ \left. \varpi^0_0(X) \theta^{\bar 0}(X)
- X(\theta^{\bar 0}(X))
\right|_{\tilde \gamma_t}  = 0\ , \quad 
\left.  \varpi^{\bar 0}_{\bar 0}(X)\theta^{0}(X)
- X(\theta^{0}(X))\right|_{\tilde \gamma_t} = 0\ .\tag6.10$$
Since $\varpi^0_0 + \varpi^{\bar 0}_{\bar 0} = 0$,
multiplying the first and the second equation 
by $\theta^{0}(X)$ and $\theta^{\bar 0}(X)$, respectively, and then
adding them together we get: 
$$\frac{d(\theta^0(X)\theta^{\bar 0}(X))}{dt} = \frac{d F^2(\dot \gamma)}{dt} 
  = 0\ .\tag 6.11$$
This shows that  the  equations (6.10) are equivalent to the 
 equations (6.5) and it concludes the proof.\qed
\enddemo
In analogy with what happens in Riemannian geometry,
it is not hard to realize that 
 a regular curve is critical 
w.r.t. the length functional if and only if, up to a parameterization,
it is a geodesic.
\bigskip
\subsubhead 6.2 Complex geodesics,
E-manifolds and Lempert manifolds
\endsubsubhead
\medskip
Let $(N,J_N)$  be a 
 complex manifold of dimension $\dim_\C N \leq \dim M = n$
and let $\imath: N\to M$ be 
a holomorphic embedding.  In analogy with the Riemannian
and Hermitian settings,  an holomorphic embedding  is said to be
{\it totally geodesic\/} whenever any geodesic $\gamma: [a,b] \to N$
of the induced  Finsler metric $F_N = F\circ \imath_*$ is embedded as  
a geodesic $\gamma' = \imath\circ \gamma$ of $M$.\par 
We give  here the concepts of complex geodesics and complex pre-geodesics. Note
that our definition of  complex geodesics coincides with that of 
segments of complex geodesics given in [AP]. They are strongly related (but different)
with the complex geodesics as defined by Vesentini in [Ve]
(see remarks in [AP], p. 129).
\medskip
\definition{Definition 6.3}
A {\it complex pre-geodesic\/} of a complex Finsler manifold $(M,J,F)$
is a totally geodesic holomorphic  
embedding $\imath: \Gamma \to M$ of a simply connected
complex curve  $(\Gamma, J_o)$. \par
A {\it complex geodesic\/} is a complex pre-geodesic
$\imath: \Gamma \to M$ such that the 
 K\"ahler metric induced on $\Gamma$ by $M$
has constant holomorphic
sectional curvature.\par
A complex Finsler manifold $(M,J,F)$ is called
{\it E-manifold\/} if
\roster
\item"i)" for any $x\in M$ and any 
vector $v\in T_xM$ there exists a 
complex geodesic $\imath: \Gamma \to M$ passing through $x$
and tangent to $v$;
\item"ii)" all complex geodesics have the same  holomorphic 
sectional curvature.
\endroster
\enddefinition
\medskip
As we mentioned in the 
Introduction, the examples of E-manifolds we are mainly interested in
are the Lempert
manifolds (see Definition 1.1): they are  E-manifolds with
complex geodesics
of   holomorphic sectional curvature   
$-4$.  Other interesting examples of complete E-manifolds, 
with complex geodesics with non-negative holomorphic sectional curvature, 
are given by the classification of non-negatively curved Kahler-Finsler 
manifolds given by Abate and Patrizio in [AP1]. \par
\medskip
The goal of this subsection is to give some properties 
on the torsion and the curvature, which gives
a complete characterization of the 
E-manifolds. In the following Theorem 6.7, we will show that the E-manifolds
are exactly the 
 complex Finsler manifolds, which are
 {\it geodetically torsion-free
and with  constant
holomorphic sectional curvature\/} (see
Definition 6.4 below).\par 
Notice  that what we call
 {\it geodetically torsion-free Finsler manifolds\/} coincide
with the manifolds that  Abate and Patrizio
christened {\it weakly-K\"ahler Finsler manifolds\/} 
(see [AP]). \par
\medskip
We begin with the equations of complex pre-geodesics
and the complex geodesics of a complex Finsler manifold.\par
\medskip
Let $\imath: \Gamma \to M$ be a holomorphic embedding of a complex
curve $\Gamma$ and let $F_\Gamma = F\circ \imath_*$ the 
induced Finsler metric on $\Gamma$. 
We call {\it adapted unitary frame of $(\Gamma, \imath)$\/} any frame 
$u = \{f_0, \dots, f_{2n-1}\}\in \pi^{-1}(\imath(\Gamma)) \subset \UF$ with 
$f_0$ tangent to $\imath(\Gamma)$.
We denote by $U_{\imath}(\Gamma)$ the bundle of 
adapted unitary frames. It is immediate to realize 
that $U_{\imath}(\Gamma)/U_{n-1} = 
U_{F_\Gamma}(\Gamma)$.\par
Let us denote by $\tilde \imath: U_{\imath}(\Gamma) \to \UF$ the natural
immersion map. Then we have the following commutative diagram.
$$
\CD
 U_{\imath}(\Gamma) @>{\tilde \imath}>>
\UF  \\
@V{\hat \pi}VV @VV{\hat\pi}V\\
U_{F_\Gamma}(\Gamma) = U_{\imath}(\Gamma)/U_{n-1}@>>> 
SM = \UF/U_{n-1}\\
@V{\pi'}VV @VV{\pi'}V\\
\Gamma @>{\imath}>> M
\endCD
\tag6.14$$
Let us also define on $U_{\imath}(\Gamma)$ the following 1-forms
$$\vartheta^\alpha = \tilde \imath^* \theta^\alpha\ ,\quad
\pi^\alpha_\beta = \tilde \imath^* \varpi^\alpha_\beta\ .\tag6.15$$
If we denote by 
 $\theta_\Gamma^0$ and $\theta_\Gamma^{\bar 0}$  the
holomorphic  components of the tautological 
1-form of $U_{F_\Gamma}(\Gamma)$, then
$$\vartheta^0 \equiv \hat \pi^*\theta_\Gamma^0\ ,\quad \vartheta^\lambda \equiv 0
\ .\tag6.16$$
\proclaim{Lemma 6.4} Let $\imath: \Gamma\to M$ be an holomorphic embedding  of a 
complex curve $\Gamma$ in $(M,J,F)$ and let $\underset{\Gamma}\to
\varpi^{\alpha}_\beta$
the components of the Hermitian connection of the induced metric 
$F_\Gamma = F\circ \imath_*$ on $\Gamma$.
The embedding is totally geodesic
if and only if for any $\lambda = 1, \dots, n-1$:
\roster
\item $\pi^0_0 = \hat\pi{}^*
\underset{\Gamma}\to\varpi^0_0$;
\item $\pi^\lambda _0\equiv 0$ and 
$\pi^0_\lambda \equiv 0$;
\item $T^0_{\lambda 0}\circ \tilde \imath \equiv 0$.
\endroster
In particular if $\imath$
is a totally geodesic
holomorphic embedding, then for any $u\in  U_{F_\Gamma}(\Gamma)$ and
any $\tilde u$  of $\hat \pi^{-1}(u)
\in U_{\imath}(\Gamma)$, 
$R^0_{\lambda 0 \bar 0}|_{\tilde \imath(u)} = 0$ and
the  holomorphic sectional curvature $c|_u$ of $F_\Gamma$
is equal to 
$$c|_u = R^0_{0 0 \bar 0}|_{\tilde \imath(u)}\ .$$
\endproclaim
\demo{Proof} Let $\gamma: [a, b] \to \Gamma$ be a regular curve in 
$\Gamma$ and let $\tilde \gamma: [a,b] \to U_{\imath}(\Gamma) \subset 
\UF$ be a lift of the curve $\gamma_\imath =
\imath\circ \gamma$. By (6.16) and 
Theorem 6.2, $\gamma$ is a geodesic for the induced metric if and only 
if
$$\frac{d F_\Gamma(\dot \gamma)}{dt} = 
\frac{d F(\dot{\gamma}_\imath)}{dt} = 0\ ,
\quad
\vartheta^{\bar 0}
(\dot{\tilde \gamma})
\left(\hat \pi^*\underset{\Gamma}\to\varpi
^0_0\right)(\dot{\tilde \gamma_t}) - \left.\frac{d \vartheta^{\bar 0}
(\dot{\tilde \gamma_s})}{ds}\right|_t  = 0\ .\tag6.17$$
Using again Theorem 6.2, $\gamma_\imath$ is  a geodesic
for the Finsler metric of 
$M$  if and only if $\tilde \gamma$ verifies also
$$
\vartheta^{\bar 0}
(\dot{\tilde \gamma}_t)
\pi^0_0(\dot{\tilde \gamma_t}) - \left.\frac{d \vartheta^{\bar 0}
(\dot{\tilde \gamma}_s)}{ds}\right|_t = 0\ ,\tag6.18$$
$$\pi^0_\lambda(\dot{\tilde \gamma_t}) + T^0_{\lambda 0}|_{\tilde \gamma_t} 
\vartheta^0
(\dot{\tilde \gamma}) = 0\ ,\quad
\pi^{\bar 0}_{\bar \lambda}
(\dot{\tilde \gamma_t}) + 
\overline{T^{0}_{\lambda 0}}|_{\tilde \gamma_t} 
\vartheta^{\bar 0}
(\dot{\tilde \gamma}) = 0\ .\tag6.19$$
Therefore the embedding is totally geodesic  if and only if  (1) holds and
$$\pi^0_\lambda|_u = - T^0_{\lambda 0}|_u \vartheta^0_u\ ,\qquad
\pi^{\bar 0}_{\bar \lambda}|_u = 
- \overline{T^{0}_{\lambda 0}}|_u \vartheta^{\bar 0}_u \tag6.20$$
for any $u\in U_{\imath}(\Gamma)$. On the other hand, by the 
structure equations (5.29)  and (5.30) and by (6.16)
$$0 = d\vartheta^\lambda 
= - \pi^{\lambda}_0\wedge \vartheta^0 + \imath^* \Theta^\lambda 
= \pi^{\bar 0}_{\bar \lambda}\wedge \vartheta^0
\ .\tag6.21$$
From (6.20) and (6.21) it follows
that 
$$(\overline{T^{0}_{\lambda 
0}}\circ {\tilde \imath}) \vartheta^{\bar 0}\wedge \vartheta^0
\equiv 0$$
which implies that $\overline{T^{0}_{\lambda 0}}
\circ {\tilde \imath}\equiv 0$ since
$\vartheta^0\wedge \vartheta^{\bar 0} = \hat\pi^*(\theta_\Gamma^0
\wedge \theta^{\bar 0}_\Gamma) \neq 0$. From this, (2) and (3) are immediate.\par
The last claims follows from (1), (2) and  the 
structure equations of $U_{F_\Gamma}(\Gamma)$.\qed
\enddemo
\bigskip
We can now give the characterization of E-manifolds. Let us first
introduce some terminology.\par
\medskip
\definition{Definition 6.5} We say that 
a complex Finsler manifold $(M,J,F)$ is 
called {\it geodetically torsion-free\/} if  the 2-form
$\Theta^0$  is of the following form
$$\Theta^0 = \frac{1}{2}
T^0_{\lambda \mu} \theta^\lambda\wedge\theta^\mu\tag6.22$$
(i.e. the complex functions 
 $T^{0}_{\lambda 0}$  vanish identically).\par
$(M,J,F)$ is called  {\it with constant holomorphic sectional curvature\/} if there
exists a constant $c$ so that 
the  2-forms $\Omega^0_0$ and $\Omega^0_\lambda$ are of the form
$$\Omega^0_0 = c\theta^0\wedge \theta^{\bar 0} + 
R^0_{0\rho \bar \sigma}\theta^\rho \wedge\theta^{\bar \sigma} +
R^0_{0\rho \bar 0}\theta^\rho \wedge\theta^{\bar 0} +
R^0_{0 0 \bar \sigma}\theta^0 \wedge\theta^{\bar \sigma}\ ,$$
$$ \Omega^\lambda_0 = 
R^\lambda_{0 \rho \bar \sigma}\theta^\rho \wedge\theta^{\bar \sigma}
+ R^\lambda_{0 \rho \bar 0}\theta^\rho \wedge\theta^{\bar 0}
+ R^\lambda_{0 0 \bar \sigma}\theta^0 \wedge\theta^{\bar \sigma}\ ,
$$
$$ \Omega^0_\lambda = 
R^0_{\lambda \rho \bar \sigma}\theta^\rho \wedge\theta^{\bar \sigma}
+ R^0_{\lambda \rho \bar 0}\theta^\rho \wedge\theta^{\bar 0}
+ R^0_{\lambda 0 \bar \sigma}\theta^0 \wedge\theta^{\bar \sigma}
\tag6.23$$
(i.e. $R^0_{0 0\bar 0} \equiv c$ and $R^\lambda_{0 0\bar 0}
\equiv R^0_{\lambda 0\bar 0} \equiv 0$). \par
If $M$ has constant holomorphic sectional curvature, 
the constant  $c$ is called {\it the holomorphic
sectional curvature of $M$\/}.
\enddefinition
\medskip
\remark{Remark 6.6}  Assume that  $F$
is associated with an Hermitian metric $g$. In this case, using the 
fact that the functions $T^\alpha_{\beta \gamma}(u)$ depends
on the frame $u$ as the components of a tensor of type $(1,2)$,
it can be inferred that $F$ is geodetically torsion-free if and only if 
$g$ is torsion free and hence K\"ahler. With the same 
arguments,  it can be shown that $F$ is  
of constant holomorphic sectional curvature if and only if the Hermitian
metric $g$ is of constant holomorphic sectional curvature.
\endremark
\medskip
Here is the characterization we were looking for.
\bigskip
\proclaim{Theorem 6.7} Let  $(M,J,F)$ be  a complex Finsler manifold. 
\roster
\item"i)"
There exists  a complex pre-geodesic  through any point $x\in M$ and 
tangent to any vector $v\in T_xM$ if and only $M$ 
is geodetically torsion-free. 
\item"ii)" $(M, J, F)$ is an E-manifold if and only if it is 
geodetically torsion-free and with constant holomorphic sectional curvature.
\endroster
\endproclaim
\demo{Proof} (i) The necessity follows immediately from Lemma 6.5. 
Suppose now that $(M, J, F)$ is geodetically torsion-free and consider 
the distribution $\Cal C$ on $\UF$ given by all vectors $X\in T\UF$
such that
$$\theta^\lambda(X) = 0\ ,\quad \theta^{\bar \lambda}(X) = 0
\ ,\quad\varpi^{\lambda}_0(X) = 0\ ,
\quad \varpi^0_{\lambda}(X) = 0\tag6.24$$
for  $\lambda = 1, \dots, n-1$. Using the structure equations, 
one can check that the equations (6.24) define an integrable
distribution whose integral leaves of maximal
dimension project onto holomorphic
curves in $(M,J)$. Moreover, if $S\subset \UF$ is an integral leaf
of $\Cal C$ with corresponding
holomorphic curve
 $\Gamma = \pi(S)\subset M$  and if $\imath: \Gamma \to M$  is the
standard immersion of $\Gamma$, then $S$ is equal to the adapted
frame bundle $U_{\imath}(\Gamma)$ and
 the immersion
$\imath: \Gamma \to M$ is a totally geodesic isometric embedding.
Since there exists an integral leaf of $\Cal C$ for any frame $u\in \UF$,
this concludes the proof of (i).\par
The proof of (ii) is analogous.\qed
\enddemo
\bigskip
\remark{Remark 6.8} Equivalent characterizations 
of the E-manifolds can be also found in [Fa], [Pa] and [AP].
\endremark
\bigskip
\subsubhead 6.3 The torsion and curvature  of an E-manifold
\endsubsubhead
\bigskip
In the following last Theorem 6.7, we prove that 
the torsion  and the curvature of an E-manifold
are uniquely determined by the Finsler torsion and 
the Finsler curvatures. This implies that in order
to have a complete set of invariants for an E-manifold, 
it suffices to consider the structure functions
described in Corollary 5.7 i), ii) and iii). We also give the explicit
formulae for some components of the torsion and the curvature and 
an application of these formulae, which gives a short proof
of an Abate and Patrizio's result on K\"ahler-Finsler manifolds
with positive sectional curvature (Theorem 1.1 in [AP1]).\par
For the notation and the indexing conventions,
see \S 5.1 and \S 5.2.\par
\medskip
\proclaim{Theorem 6.9} Let $(M,J,F)$ be an E-manifold with constant 
holomorphic sectional curvature $c$. Then:
\roster
\item the torsion and the
curvature of the non-linear Hermitian connection of $M$
are uniquely determined by the structure functions $h_{\lambda \mu}$, $H_{\lambda
\mu \bar \nu}$, 
$P^\lambda_{0 \bar \mu \gamma} = \hat e_\gamma(h_{\bar \lambda \bar \mu})$
and $Q^\rho_{\sigma\lambda\bar \mu} =   H_{\bar \mu \bar \rho \sigma \lambda} 
- h_{\bar \mu \bar \rho} h_{\sigma \lambda} -
 H_{\nu \bar \mu \bar \rho }
H_{\sigma \lambda \bar \nu}$ and their first
order derivatives; in particular, 
$$R^0_{0 0 \bar 0} = c\ ,\quad
R^0_{\lambda 0 \bar 0} = R^\lambda_{0 0 \bar 0} = R^0_{0 \lambda \bar 0} = 
R^0_{0 0 \bar \lambda} = 0\ ,\quad
R^0_{\lambda \mu \bar 0} = c h_{\lambda \mu} \ ,\quad 
R^\lambda_{0 0 \bar \mu} = c h_{\bar \lambda \bar \mu} \ ,$$
$$
R^0_{\lambda 0 \bar \mu} = \frac{c}{2} (\delta_{\lambda\mu} +
h_{\lambda \rho} h_{\bar \rho \bar \mu})\ ,\quad 
R^\lambda_{0 \mu \bar 0} = \frac{c}{2} (\delta_{\lambda\mu} +
h_{\bar \lambda \bar \rho} h_{\rho \mu})\ ,\quad
R^0_{0\lambda \bar \mu} = \frac{c}{2} (\delta_{\lambda\mu} - 
h_{\mu \rho} h_{\bar \rho \bar \lambda})\ ,$$
$$R^\lambda_{\mu 0 \bar 0} =  
\frac{c}{2}(\delta_{\lambda \mu} -
h_{\nu \mu}  h_{\bar \nu \bar \lambda}) - \hat e_{\bar 0}(h_{\mu \nu})
\hat e_{0}(h_{\bar \lambda \bar \nu})\ .$$
\item if $c \neq 0$,  the $0$-th component of the 
torsion $\Theta^0 = T^0_{\beta \gamma} \theta^\beta\wedge \theta^\gamma$
  vanishes identically and the whole set of components of the torsion
is given by the following expressions:
$$T^0_{\beta \gamma} = 0\ ,\quad T^\beta_{0 \gamma} = 
- \hat e_0(h_{\bar \beta \bar \nu}) h_{\nu \gamma}\ ,
\quad  T^\alpha_{\beta\gamma} = 
 \hat e_\gamma(h_{\bar \alpha \bar \nu}) h_{\nu \beta}
- \hat e_\beta(h_{\bar \alpha \bar \nu}) h_{\nu \gamma}\ .$$
\item if $c>0$ and the functions $\hat e_0(T^\lambda_{0 \lambda})$ vanish identically 
for any 
$\lambda = 1, \dots, n-1$, 
then  $F$ is associated with a K\"ahler metric of 
constant holomorphic sectional curvature; in particular, if $(M,J,F)$ 
 is also simply connected
and complete, then it is biholomorphic to $\C P^n$.
\endroster
\endproclaim
\demo{Proof}  (1) The proof is based on iterated
use of the identities (5.36), the Bianchi identities and the 
Jacobi identities applied to three vector fields $v_1$, $v_2$ and $v_3$ on $\UF$, where 
$v_1$ and $v_2$ are of the form $\hat e_\alpha$ or $\hat e_{\bar \alpha}$ and
 $v_3$ is a vector field of the form $\tilde e_\lambda$ or $\tilde e_{\bar \lambda}$.
The arguments are simple and straightforward and we are going to 
show only how to determine the expressions for the components 
$R^\alpha_{\beta \gamma \bar \delta}$ were at least two indices are equal
to $0$ or $\bar 0$. The way to determine all other components of the curvature 
and of the components of the torsion are analogous.\par
By hypotheses, for any $\lambda = 1, \dots n-1$,
 $R^0_{00\bar 0} =  c$, $R^0_{\lambda 0\bar 0} =  R^\lambda_{00\bar 0} =  0$
and $T^0_{0\lambda} = 0$. Then, 
from the Bianchi identities (5.42) we  get
$$R^0_{0 \lambda \bar 0} - R^0_{\lambda 0\bar 0}  -
\hat e_{\bar 0}(T^0_{0 \lambda}) - H_{\bar 0 \mu 0} R^\mu_{0\lambda \bar 0}
+ H_{\bar 0 \mu \lambda} R^\mu_{0 0 \bar 0} = 
R^0_{0 \lambda \bar 0} - R^0_{\lambda 0\bar 0} = 0\ .\tag6.25$$
On the other hand, by (5.36)
$$R^\alpha_{\beta \gamma \bar \delta} = \overline{R^\beta_{\alpha \delta \bar \gamma}}
\tag6.26$$
From (6.25) and (6.26), we conclude that $
R^0_{0 0 \bar \lambda} = R^0_{0 \lambda \bar 0} = 0$.\par
Now, using the notation of \S 5, by the Jacobi identities
we have  
$$\omega^0_0([\tilde e_\lambda, [\hat e_\mu, \hat e_{\bar 0}]])
+ \omega^0_0([\hat e_\mu, [\hat e_{\bar 0}, \tilde e_{\lambda}]]) + 
\omega^0_0([\hat e_{\bar 0}, [\tilde e_{\lambda}, \hat e_{\mu}]]) = $$
$$ = R^0_{\lambda \mu \bar 0} - \tilde e_\lambda(R^0_{0 \mu \bar 0})
- H_{\alpha \lambda \mu} R^0_{0 \alpha\bar 0} = 
 R^0_{\lambda \mu \bar 0} - c h_{\lambda \mu} = 0\tag 6.27$$
From (6.27) and (6.26), it follows also that 
$R^\lambda_{0 0 \bar \mu}  = \overline{R^0_{\lambda \mu \bar 0}} = c h_{\bar \lambda
\bar \mu}$. \par
Let us use again the first Bianchi identities and the Jacobi identities:
$$R^0_{\lambda 0 \bar \mu} - R^0_{0 \lambda \bar \mu}
- \hat e_{\bar \mu}(T^0_{\lambda 0}) - H_{\bar 0 \nu \lambda} R^{\nu}_{0 0 \bar \mu}
+  H_{\bar 0 \nu 0} R^{\nu}_{0 \lambda \bar \mu} = 
R^0_{\lambda 0 \bar \mu} - R^0_{0 \lambda \bar \mu}
- c h_{\nu \lambda}  h_{\bar \nu
\bar \mu} = 0\tag6.28$$
$$\omega^0_0([\tilde e_\lambda,[\hat e_0, \hat e_{\bar \mu}]])  + 
\omega^0_0([\hat e_0,[\hat e_{\bar \mu}, \tilde e_{\lambda}]]) + 
\omega^0_0([\hat e_{\bar \mu},[\tilde e_{\lambda}, \hat e_{0}]]) = $$
$$ - \tilde e_{\lambda}(R^0_{0 0 \bar \mu}) + R^0_{\lambda 0 \bar \mu}
- R^0_{0 0 \bar 0}\delta_{\lambda \mu} + R^0_{0 0 \lambda \bar \mu} = 
R^0_{\lambda 0 \bar \mu}
+ R^0_{0 \lambda \bar \mu} - c\delta_{\lambda \mu} = 0\tag6.29$$
From (6.28), (6.29) and (6.26), it follows that 
$$R^0_{\lambda 0 \bar \mu} = \overline{R^\lambda_{0  \mu \bar 0}} = 
\frac{c}{2}(\delta_{\lambda \mu} + 
h_{\nu \lambda}  h_{\bar \nu \bar \mu})\ , \qquad 
R^0_{0 \lambda \bar \mu} = 
\frac{c}{2}(\delta_{\lambda \mu} -
h_{\nu \lambda}  h_{\bar \nu \bar \mu})\ .$$
Using again the Jacobi identities,
$$\omega^{\bar \mu}_{\bar O}([\tilde e_{\bar \lambda}, [\hat e_0, 
\hat e_{\bar 0}]]) + 
\omega^{\bar \mu}_{\bar O}([\hat e_{0}, [\hat e_{\bar 0}, 
\tilde e_{\bar \lambda}]]) + 
\omega^{\bar \mu}_{\bar O}([\hat e_{\bar 0}, [\tilde e_{\bar \lambda}, 
\hat e_{0}]]) = $$
$$ R^\lambda_{\mu 0 \bar 0} - c \delta_{\lambda\mu} 
 + \frac{c}{2}(\delta_{\lambda \mu} + 
h_{\nu \mu}  h_{\bar \nu \bar \lambda}) + \hat e_{\bar 0}(h_{\mu \nu})
\hat e_{0}(h_{\bar \lambda \bar \nu}) = 0\ .$$
\medskip
(2) Assume $c \neq 0$. From the Bianchi identity (5.43) and (1),
$$\hat e_{\beta}(R^0_{0 \gamma \bar 0}) - \hat e_{\gamma}(R^0_{0 \beta \bar 0}) 
+ R^0_{0 \delta \bar 0} T^\delta_{\beta \gamma} = c T^0_{\beta \gamma} = 0\ ,$$
and this implies $T^0_{\beta \gamma} = 0$. Then using again the Jacobi identities, by the 
vanishing of $\sum_{cyclic\ perm.}\theta^0([v_i, [v_j, v_k]]) = 0$
when $v_1 = \hat e_\alpha$, $v_2 = \hat e_\nu$ and $v_3 = \tilde e_{\bar \lambda}$,
one obtains the remaining expressions for $T^\lambda_{0 \nu}$ and 
for $T^\lambda_{\mu\nu}$.\par
\medskip 
(3) From the Bianchi identity (5.36), we have that
$$R^\lambda_{0 \lambda \bar 0} - R^\lambda_{\lambda 0 \bar 0} - 
\hat e_{\bar 0}(T^\lambda_{0 \lambda}) = 0\ .$$
By the expressions for the curvature components given in (1),  this becomes
$$\frac{c}{2}(1 + \sum_{\rho = 1}^{n-1} |h_{\lambda \rho}|^2) 
- \frac{c}{2}(1 - \sum_{\rho = 1}^{n-1} |h_{\lambda \rho}|^2) + 
\sum_{\rho = 1}^{n-1} |\hat e_{\bar 0}(h_{\lambda \rho})|^2  - 
\hat e_{\bar 0}(T^\lambda_{0 \lambda})  =$$
$$ = 
\sum_{\rho = 1}^{n-1} \left(c|h_{\lambda \rho}|^2 + 
|\hat e_{\bar 0}(h_{\lambda \rho})|^2\right) 
- \hat e_{\bar 0}(T^\lambda_{0 \lambda}) = 0\ .$$
This implies that, if $c>0$ and 
$\hat e_{\bar 0}(T^\lambda_{0 \lambda}) = 0$ for any $\lambda$, then 
$\hat e_{\bar 0}(h_{\lambda \rho}) = 
h_{\lambda \rho} = 0$ for any $\lambda$ and $\rho$. Therefore, 
by Lemma 2.5, $F$
is associated with an Hermitian metric $g$, which is geodetically 
torsion free and with constant holomorphic sectional curvature. By Remark 6.6,
we obtain that $(M, J, g)$ is K\"ahler and locally isometric to $\C P^n$. 
The conclusion  follows from  
standard facts on  complex space forms.\qed 
\enddemo
\bigskip
\bigskip
\head Appendix
\endhead 
\medskip
We recall here the Cartan-Sternberg theorem on the
local automorphisms of an absolute parallelism.  The theorem was first 
proved for real analytic vector fields by E. Cartan and in this case it is a 
corollary 
of  Cartan-K\"ahler theorem (see e.g. [BCG]). Later it was proved
 by S. Sternberg for 
smooth vector fields ([St]).  \par
\smallskip
Before stating the theorem  we need some preliminaries.\par
Let $\sigma  = \{X_1,  \dots, X_n\}$ be an absolute parallelism on a 
manifold $N$.
The {\it structure functions\/} of $\sigma$ are the smooth functions
$c^i_{jk}$ defined by 
$$[X_j, X_k] = \sum_{i=1}^n c^i_{jk} X_i\ .$$
Let us also denote by $c^i_{jk, m_1 \dots m_r}$ the smooth functions defined 
inductively on $r$ as
$$c^i_{jk, m_1} = X_{m_1}(c^i_{jk})\ ,\qquad c^i_{jk, m_1 \dots m_r} = 
X_{m_r}(c^i_{jk,m_1 \dots m_{r-1}})\ .$$
Finally, for any integer $\alpha >0$ let $\F\alpha$ be
the family of smooth functions
$$\F\alpha = \{ c^i_{jk}, c^i_{jk,m_1}, \dots, 
c^i_{jk, m_1 \dots m_\alpha}\}$$
and call $Q_\alpha$ the number of functions in the set $\F\alpha$. 
We  consider
$\F\alpha$ as the set of components of a smooth map from $M$
into $\R^{Q_\alpha}$. A point $p\in M$ is called a {\it regular point
for $\sigma$\/} if there exists two integers $s$ and $r$ such that 
$\rank \F s = \rank \F{s+1} = r$ 
at all points
of a neighborhood $U_p$ of $p$. \par
If $s$ is the smallest integer
such that this occurs, then $s$ and $r$ are called {\it order\/}  
and {\it rank\/} of the regular point $p$, respectively. \par
It can be shown that 
 $\rank \F \alpha = r$ for all $\alpha \geq
s$ and that
there exists a system of coordinates $\{ x_1, \dots, x_n\}
: U_p \to \R^n$ such that all maps $\F \alpha|_{U_p}$, $\alpha >0$,  depend only 
on the first $r$ coordinates $\{ x_1, \dots, x_r\}$ (see [St]).
 Such a system of coordinates
is called {\it adapted to the absolute parallelism\/}. \par
For any 
the  adapted system of coordinates $\{ x_1, \dots, x_n\}$ on a neighborhood
$U_p$, we call  {\it  slice of $U_p$\/}
any set of the form
$$S_{(c_1, \dots, c_r)} = \{ q \in U\ :\ 
x_1(q) = c_1, \dots , x_r(q) = c_r\ \}\ ,$$
for some $(c_1, \dots, c_r) \in \R^r$.
\medskip
\proclaim{Theorem A1} (Cartan -Sternberg) Let 
$\sigma = 
\{ X_1, \dots, X_n\}$ be an absolute parallelism on $M$ and let 
$p, p'\in M$ be two regular point of ranks $r_p$ and $r_{p'}$ 
and orders $s_p$ and $s_{p'}$, respectively. Let also
 $U$ and $U'$ be two neighborhoods
of $p$ and $p'$, respectively, which admit two adapted systems of 
coordinates $\{x_i\}$ and $\{x'_i\}$.\par
If $q$ is a point of the slice  $S_{(c_1, \dots, c_r)}\subset  U$ and $q'$
is a point of the slice $S_{(c'_1, \dots, c'_r)}\subset U'$,
 there exists
a local diffeomorphism $f: U\to U'$ such that $f(q) = q'$ 
and $f_*(X_i) = X_i$ for all $i = 1, \dots, n$, if and only if $r_p = r_{p'} = r$, 
$s_p = s_{p'} = s$ and 
$\F s|_{S_{(c_1, \dots, c_r)}} \equiv 
\F s|_{S_{(c'_1, \dots, c'_r)}}$.\par
In particular, if $U = U'$,  there exists
a local diffeomorphism $f: U\to U$ such that $f(q) = q'$ 
and $f_*(X_i) = X_i$ for all $i = 1, \dots, n$ if and only if 
$q$ and $q'$ belong to the same slice $S_{(c_1, \dots, c_r)}$ for 
some $(c_1, \dots, c_r) \in \R^r$.
\endproclaim

\enddocument
\bye